\documentclass[10pt,twoside,a4paper,reqno]{amsart}
\usepackage{amsfonts,amsthm,latexsym,amsmath,amssymb,amscd,epsfig}
\usepackage{graphics,graphicx}

\theoremstyle{plain}
\newtheorem{theo+}           {Theorem}
\newtheorem{prop+}           {Proposition}
\newtheorem{coro+}           {Corollary}
\newtheorem{lemm+}           {Lemma}
\newtheorem{conjecture}      {Conjecture}
\newtheorem*{nutt}           {Nuttall's Conjecture}
\newtheorem*{pade}           {Pad\'e (Baker-Gammel-Wills) Conjecture}

\theoremstyle{definition}
\newtheorem{defi+}           {Definition}
\newtheorem{problem}         {Problem}

\newtheorem*{pb-prime}       {Problem 1$\mathbb{'}$}
\newtheorem*{ack}            {Acknowledgements}

\newtheorem{not+}            {Notation}

\theoremstyle{remark}
\newtheorem{rema+}           {Remark}
\newtheorem*{expl1}          {Explanations to Figure~\ref{fig1}}
\newtheorem*{expl2}          {Explanations to Figure~\ref{fig2}}
\newtheorem{example}         {Example}

\newenvironment{theorem}{\begin{theo+}}{\end{theo+}}
\newenvironment{proposition}{\begin{prop+}}{\end{prop+}}
\newenvironment{corollary}{\begin{coro+}}{\end{coro+}}
\newenvironment{lemma}{\begin{lemm+}}{\end{lemm+}}
\newenvironment{remark}{\begin{rema+}}{\end{rema+}}
\newenvironment{definition}{\begin{defi+}}{\end{defi+}}
\newenvironment{notation}{\begin{not+}}{\end{not+}}

\newcommand{\al}{\alpha}
\newcommand{\be}{\beta}
\newcommand{\ga}{\gamma}
\newcommand{\Ga}{\Gamma}
\newcommand{\ze}{\zeta}
\newcommand{\la}{\lambda}

\newcommand{\Si}{\Sigma}
\newcommand{\bC}{\mathbb C}
\newcommand{\bCP}{\mathbb {CP}}
\newcommand{\bR}{\mathbb R}
\newcommand{\bN}{\mathbb N}
\newcommand{\bZ}{\mathbb Z}
\newcommand{\bD}{\mathbb D}
\newcommand{\A}{\mathcal A}
\newcommand{\D}{\mathcal D}
\newcommand{\C}{\mathcal C}

\newcommand{\N}{\mathcal N}
\newcommand{\cS}{\mathcal S}
\newcommand{\SG}{\mathcal{SG}}
\newcommand{\eps}{\epsilon}
\newcommand {\cP} {\mathcal P}
\newcommand {\Up} {\Upsilon}
\newcommand{\Om}{\Omega}
\newcommand{\bx}{\mathbf x}
\newcommand{\bu}{\mathbf u}
\newcommand{\bv}{\mathbf v}

\newcommand{\bq}{\mathbf q}

\newcommand{\cV}{\mathcal V}

\newcommand{\De}{\Delta}

\numberwithin{equation}{section}

\begin{document}

\title[Rational approximation of algebraic functions]
{On rational approximation \\ of algebraic functions}
\author[J.~Borcea]{Julius Borcea$^*$}
\thanks{$^*$ Corresponding author}
\address{Department of Mathematics, Stockholm University, SE-106 91 Stockholm,
      Sweden}
\email{julius@math.su.se}
\author[R.~B\o gvad]{Rikard B\o gvad}
\address{Department of Mathematics, Stockholm University, SE-106 91 Stockholm,
      Sweden}
\email{rikard@math.su.se}
\author[B.~Shapiro]{Boris Shapiro}
\address{Department of Mathematics, Stockholm University, SE-106 91 Stockholm,
       Sweden}
\email{shapiro@math.su.se}

\begin{abstract} We construct a new scheme of approximation of any
multivalued algebraic function $f(z)$ by a sequence
$\{r_{n}(z)\}_{n\in \bN}$ of rational functions. The latter sequence is
generated by a recurrence relation which is completely determined by the
algebraic equation satisfied by $f(z)$. Compared to the usual Pad\'e
approximation our scheme has a number of advantages, such as
simple computational procedures
that allow us to prove natural analogs of the Pad\'e Conjecture and
Nuttall's Conjecture for the sequence
$\{r_{n}(z)\}_{n\in \bN}$ in the complement
$\bCP^1\setminus \D_{f}$, where $\D_{f}$ is the union of a finite
number of segments of real algebraic curves and finitely many isolated points.
In particular, our construction makes it possible to control the behavior of
spurious poles and to describe the
asymptotic ratio distribution of the family $\{r_{n}(z)\}_{n\in \bN}$. As an
application we settle the so-called $3$-conjecture of Egecioglu {\em et al}
dealing with a $4$-term recursion related to a polynomial Riemann Hypothesis.
\end{abstract}

\maketitle

\section{Introduction and main results}\label{s1}

Rational approximants of analytic functions and the asymptotic distribution
of their zeros and poles are of central interest in many
areas of mathematics and physics. For the class of algebraic functions these
questions have attracted special attention due to their
important applications ranging from the convergence theory of Pad\'e
approximants \cite{St1}--\cite{St6} and the theory of general
orthogonal polynomials \cite{Ch,Ne,ST,Sz} to
statistical mechanics \cite{So1,So2}, complex Sturm-Liouville
problems \cite{BBS}, inverse scattering theory and quantum field
theory \cite{BG-M}.

The main purpose of this paper is to give a simple and direct construction of
families of rational functions converging to a certain branch of an arbitrary
(multivalued) algebraic function $f(z)$. While the usual Pad\'e approximation
requires the knowledge of the Taylor expansion of $f(z)$ at  $\infty$, our
scheme is based only on the algebraic equation satisfied by
$f(z)$ and has therefore an essentially different range of applications. 
Indeed, let
\begin{equation}\label{eq:df}
P(y,z)=\sum_{i=0}^{k}P_{k-i}(z)y^{k-i}
\end{equation}
denote the irreducible polynomial in $(y,z)$ defining $f(z)$, that is,
$P(f(z),z)=0$. Note that $P(y,z)$ is uniquely defined up to a scalar factor.
Let us rewrite~\eqref{eq:df} as
\begin{equation}\label{eq:df1}
          -y^k=\sum_{i=1}^{k}\frac{P_{k-i}(z)}{P_{k}(z)}y^{k-i}
\end{equation}
          and consider the
          associated recursion of length $k+1$ with rational
          coefficients
\begin{equation}\label{eq:rec}
	-q_{n}(z)=\sum_{i=1}^{k}\frac{P_{k-i}(z)}{P_{k}(z)}q_{n-i}(z).
\end{equation}
Choosing any initial $k$-tuple  of rational functions
	  $IN=\{q_{0}(z),\ldots,q_{k-1}(z)\}$ one can generate a family
	  $\{q_{n}(z)\}_{n\in \bN}$ of rational functions
satisfying~\eqref{eq:rec} for all $n\ge k$ and coinciding with the entries of
	  $IN$ for $0\le n\le k-1$. The main object of study of this paper
is the
	  family $\{r_{n}(z)\}_{n\in \bN}$, where
	  $r_{n}(z)=\frac{q_{n}(z)}{q_{n-1}(z)}$.

In order to formulate our results we need several additional definitions.
Consider first a usual recurrence relation of length $k+1$  with
constant coefficients
\begin{equation}\label{eq:Basic}
         -u_{n}=\al_{1}u_{n-1}+\al_{2}u_{n-2}+\ldots+\al_{k}u_{n-k},
\end{equation}
where $\al_{k}\neq 0$.

\begin{definition}\label{def4}
        The {\em asymptotic symbol equation} of
recurrence~\eqref{eq:Basic} is given by
\begin{equation}\label{eq:Char}
	t^k+\al_{1}t^{k-1}+\al_{2}t^{k-2}+\ldots +\al_{k}=0.
\end{equation}
The left-hand side of the above equation is called the
{\em characteristic polynomial} of recurrence~\eqref{eq:Basic}. Denote the
roots of~\eqref{eq:Char} by
$\tau_{1},\ldots, \tau_{k}$ and call them the {\em spectral numbers} of the
recurrence.
\end{definition}

\begin{definition}\label{def5}
Recursion~\eqref{eq:Basic} and its characteristic
polynomial are said to be of {\em dominant type} or
{\em dominant} for short if there exists a
unique (simple) spectral number $\tau_{max}$ of this recurrence relation
satisfying $|\tau_{max}|=\max_{1\le i\le k}|\tau_i|$.
Otherwise~\eqref{eq:Basic} and~\eqref{eq:Char} are said to be of
{\em nondominant type} or just {\em nondominant}. The number
$\tau_{max}$ will be referred to as the
{\em dominant spectral number} of~\eqref{eq:Basic} or the {\em dominant root}
of~\eqref{eq:Char}.
\end{definition}

The following theorem may be found in~\cite[Ch.~4]{St}.

\begin{theorem}\label{th:ordrec}
Let $k\in \bN$ and consider a $k$-tuple $(\al_{1},\ldots,\al_{k})$ of
complex numbers with $\al_{k}\neq 0$. For any function
$u:\bZ_{\ge 0}\to \bC$ the following conditions are equivalent:
\begin{enumerate}
\item[(i)] $\sum_{n\ge 0}u_{n}t^n=\frac {Q_{1}(t)}{Q_{2}(t)}$,
where $Q_{2}(t)=1+\al_{1}t+\al_{2}t^2+\ldots+\al_{k}t^k$ and $Q_{1}(t)$ is a
polynomial in $t$ whose degree is smaller than $k$.
\item[(ii)] For all $n\ge k$ the numbers $u_{n}$ satisfy the
$(k+1)$-term recurrence relation given by~\eqref{eq:Basic}.
\item[(iii)] For all $n\ge 0$ one has
\begin{equation}\label{eq:leadasymp}
	 u_{n}=\sum_{i=1}^rp_{i}(n)\tau^n_{i},
\end{equation}
where $\tau_{1},\ldots,\tau_{r}$ are the distinct spectral numbers
of~\eqref{eq:Basic} with multiplicities $m_{1},\ldots,m_{r}$, respectively,
and $p_{i}(n)$ is a polynomial in the variable $n$ of degree at most
$m_{i}-1$ for $1\le i\le r$.
\end{enumerate}
\end{theorem}

Note that by Definition~\ref{def5} the dominant
spectral number $\tau_{max}$ of any dominant recurrence relation has
multiplicity one.

\begin{definition}\label{df:slow}
        An initial $k$-tuple of complex numbers
         $\{u_{0},u_{1},\ldots,u_{k-1}\}$ is called {\em fast growing} with
respect
to a given dominant recurrence of the form~\eqref{eq:Basic} if the
         coefficient $\kappa_{max}$ of $\tau_{max}^n$ in~\eqref{eq:leadasymp} is
nonvanishing, that is, $u_{n}=\kappa_{max}\tau_{max}^n+\ldots$ with
$\kappa_{max}\neq 0$. Otherwise the $k$-tuple
$\{u_{0},u_{1},\ldots,u_{k-1}\}$ is said to be {\em slow growing}.
\end{definition}

Let $\cP_{k}=\left\{a_{k}y^k+a_{k-1}y^{k-1}+\ldots +a_{0}\mid a_i\in
\bC,0\le i\le k\right\}$
denote the linear space of all polynomials of degree at most $k$ with
complex coefficients.

\begin{definition}\label{def7}
        The real hypersurface $\Xi_{k}\subset \cP_{k}$ obtained as the
        closure
         of the set of all nondominant polynomials is called the
         {\em standard equimodular discriminant}.
For any family
$$\Gamma(y,z_{1},\ldots,z_{q})
=\left\{a_{k}(z_{1},\ldots,z_{q})y^k+a_{k-1}(z_{1},\ldots,z_{q})y^{k-1}
+\ldots+a_{0}(z_{1},\ldots,z_{q})\right\}$$
of irreducible polynomials of degree at most $k$ in the variable $y$ we
define the {\em induced equimodular
        discriminant} $\Xi_{\Gamma}$
to be the set of all parameter values
        $(z_{1},\ldots,z_{q})\in \bC^q$ for which the corresponding
polynomial in $y$ is nondominant. Given an algebraic function $f(z)$
defined by~\eqref{eq:df} we denote by $\Xi_{f}$ the induced equimodular
discriminant of~\eqref{eq:df1} considered as a family
of  polynomials in the variable $y$.
\end{definition}

	  \begin{example}
For $k=2$ the equimodular discriminant $\Xi_{2}\subset
        \cP_{2}$  is the real hypersurface consisting of all solutions to
$\eps a_{1}^2-4 a_{0} a_{2}=0$,
        where $\eps$ is a real parameter with values in $[1,\infty)$. More
        information on $\Xi_{k}$  can be found in \cite{Bi} and \cite {BSh}.
\end{example}

\begin{definition}\label{def-pl}
Given an algebraic function $f(z)$ defined by~\eqref{eq:df} and an initial
$k$-tuple
of rational functions $IN=\{q_{0}(z),\ldots,q_{k-1}(z)\}$ we define
the {\em pole locus} $\Up_{f,IN}$ associated with the data $(f,IN)$ to be
the union between the zero set of the polynomial $P_{k}(z)$ and the
sets of all poles of $q_{i}(z)$ for $0\le i\le k-1$.
\end{definition}

As we explain in \S 3, one of the fascinating features of Pad\'e approximants
is the complexity of their convergence theory. The main challenges of this
theory are the Pad\'e (Baker-Gammel-Wills) Conjecture and
Nuttall's Conjecture (cf.~\cite{BG-M,St5}). The general version of the former
was recently disproved by D.~Lubinsky~\cite{Lu}. Subsequently,
V.~Buslaev~\cite{Bu} constructed counterexamples to this conjecture for some
special algebraic (hyperelliptic) functions. Nevertheless, a question of
central interest in many applications is whether the Pad\'e Conjecture could
hold for certain classes of algebraic functions. The first main result of this
paper shows that a natural
analog of the Pad\'e Conjecture is always true for our approximation scheme:

\begin{theorem}\label{th:Exist}
Let $f(z)$ be an algebraic function defined by~\eqref{eq:df}. For any
initial $k$-tuple of rational
functions $IN=\{q_{0}(z),\ldots,q_{k-1}(z)\}$ there exists a finite set
$\Si_{f,IN}\subset \bCP^1\setminus (\Xi_{f}\cup \Up_{f,IN})$
such that
$$r_{n}(z)=\frac{q_{n}(z)}{q_{n-1}(z)}\rightrightarrows y_{dom}(z)\text{ in }
\bCP^1 \setminus \D_{f}\text{ as }n\to\infty,$$
where $y_{dom}(z)$ is the dominant root of equation~\eqref{eq:df1},
$\D_{f}=\Xi_{f}\cup\Up_{f,IN}\cup \Si_{f,IN}$, and
     $\rightrightarrows$ stands for uniform convergence on
     compact subsets of $\bCP^1 \setminus \D_{f}$.
\end{theorem}

The set $\Si_{f,IN}$ consists precisely of those points $z\in \Omega$ such
that the initial $k$-tuple $IN=\{q_{0}(z),\ldots,q_{k-1}(z)\}$ is slowly
growing with respect to recurrence~\eqref{eq:rec} evaluated at $z$
(cf.~Definition~\ref{df:slow}). This motivates the following definition.

\begin{definition}\label{Sigma}
The set $\Si_{f,IN}$ is called the {\em set of slow growth} associated with
the data $(f,IN)$.
\end{definition}

Although $\Si_{f,IN}$ is the set of solutions of an implicit equation
(cf.~\eqref{impl-sg}) it turns out that the set of slow growth is actually
a subset of a certain discriminantal set $\cS$ which is given
explicitly in terms of the data $(f,IN)$ (see the proof of Lemma~\ref{descS}).

\begin{remark}
In the final stage of preparation of this paper
the authors noticed that some weaker results similar to
Theorem~\ref{th:Exist} appeared in \cite{BKW1,BKW2}. The latter deal with
pointwise limits of certain families of recursively defined polynomials and
are often
used in modern statistical physics, see e.g.~\cite {So1} and references
therein.
\end{remark}

The next result describes the rate of convergence of the sequence
$\{r_{n}(z)\}_{n\in \bN}$ to $y_{dom}(z)$, which proves in particular the
uniform convergence stated in Theorem~\ref{th:Exist}. Given
$0<\eps\ll 1$  set
$\Theta_{\epsilon}=\bCP^1\setminus \mathcal O_{\epsilon}$, where
$\mathcal O_{\epsilon}$ is the $\epsilon$-neighborhood of
$\D_{f}=\Xi_{f}\cup\Up_{f,IN}\cup \Si_{f,IN}$ in the spherical
metric on $\bCP^1$.

\begin{theorem}\label{th:rate}
For any
sufficiently small $\epsilon >0$  the rate of convergence
of $r_{n}(z)\rightrightarrows
y_{dom}(z)$ in $\Theta_{\epsilon}$ is exponential, that is, there exist
constants $\mathfrak M>0$ and $q\in (0,1)$ such that
$\vert r_{n}(z)-y_{dom}(z)\vert\le \mathfrak M q^n$ for all $z\in
\Theta_{\epsilon}$.
\end{theorem}

\begin{definition}\label{def1}
Given a  meromorphic function $g$ in some open set $\Omega\subseteq \bC$ we
construct its (complex-valued) {\em residue distribution} $\nu_{g}$ as
follows. Let $\{z_{m}\mid m\in \bN\}$ be the (finite or infinite) set of all
the poles of $g$ in $\Omega$. Assume that the Laurent
         expansion of $g$ at $z_{m}$ has the form $g(z)=\sum_{-\infty<l\le
         l_{m}}\frac{A_{m,l}}{(z-z_{m})^{l}}$.
         Then the distribution $\nu_{g}$ is given by
\begin{equation}\label{resid}
         \nu_{g}=\sum_{m\ge 1}\left(\sum_{1\le l\le
         l_{m}}\frac{(-1)^{l-1}}{(l-1)!}A_{m,l}\frac{\partial^{l-1}}{\partial
         z^{l-1}}\delta_{z_{m}}\right),
\end{equation}
         where $\delta_{z_{m}}$ is the Dirac mass at $z_{m}$ and the sum
        in the right-hand side of~\eqref{resid} is meaningful as a distribution
in $\Omega$ since it is
         locally finite in $\Omega.$
\end{definition}

\begin{remark}
The distribution $\nu_{g}$ is a complex-valued {\em measure}
if and only if $g$ has all simple poles, see \cite[p.~250]{BeG}. If the
latter holds then in the notation of Definition~\ref{def1} the value of this
complex measure at $z_{m}$ equals $A_{m,1}$, i.e., the residue of $g$ at
$z_{m}$.
\end{remark}

\begin{definition}\label{c-trans}
Given an integrable complex-valued distribution $\rho$ in $\bC$ we define its
Cauchy transform $\C_{\rho}(z)$ as
$$\C_{\rho}(z)=\int_{\bC}\frac {d\rho(\xi)}{z-\xi}.$$
\end{definition}

It is not difficult to see that the distribution $\rho$ itself may be
restored from its Cauchy
       transform $\C_{\rho}(z)$ by
\begin{equation}\label{res-dis}
\rho=\frac 1 \pi \frac {\partial \C_{\rho}}{\partial \bar z},
\end{equation}
       where $\C_{\rho}(z)$ and
       its $\frac{\partial}{\partial \bar z}$-derivative are understood as
       distributions. As is well-known, any
meromorphic
       function $g$ defined in the whole complex plane $\bC$ is the Cauchy
transform of
       its residue distribution $\nu_{g}$  if the condition
       $\int_{\bC}d\nu_{g}(\xi)<\infty$ holds, see e.g.~\cite[p.~261]{BeG}.

\begin{definition}\label{def3}
Given a  family $\{\phi_{n}(z)\}_{n\in \bN}$ of smooth functions defined in
        some open set $\Omega\subseteq \bC$ one calls the limit
$\Phi(z)=\lim_{n\to
        \infty}\frac{\phi_{n+1}(z)}{\phi_{n}(z)}$ the {\em asymptotic
ratio} of the
family, provided that this limit exists in some open subset of $\Omega$.
If $\{\phi_{n}(z)\}_{n\in \bN}$ consists of analytic functions and $\nu_n$
denotes the
residue distribution of the meromorphic function
$\frac{\phi_{n+1}(z)}{\phi_{n}(z)}$ in $\Omega$ then the limit
        $\nu=\lim_{n\to\infty}\nu_{n}$ (if it exists in the sense
        of weak convergence) is called  the
{\em asymptotic ratio distribution} of the family.
\end {definition}

Our next main result describes the support and the density of the asymptotic
ratio distribution associated with sequences of rational functions
constructed by using our approximation scheme.

\begin{theorem}\label{th:2}
Let $f(z)$ be an algebraic function defined by~\eqref{eq:df} and fix a
generic initial $k$-tuple of rational functions
$IN=\{q_{0}(z),\ldots,q_{k-1}(z)\}$. If $\nu_{n}$ denotes the residue
distribution of
$r_{n}(z)$ and $\nu=\lim_{n\to \infty}\nu_{n}$ is the asymptotic ratio
distribution of the family $\{q_{n}(z)\}_{n\in \bN}$ then the following holds:
\begin{enumerate}
\item[(i)] The support of $\nu$ does not depend on the set of slow growth 
$\Sigma_{f,IN}$. More precisely, the asymptotic ratio distribution $\nu$ 
vanishes in a sufficiently small neighborhood of each point in $\Sigma_{f,IN}$.
\item[(ii)] Suppose that there exists a nonisolated
point $z_0\in\Xi_{f}$ such that equation \eqref{eq:df1} considered at
$z_0$ has the property that among its roots with maximal
absolute value there are at least two with the same maximal multiplicity
(compare with Definition~\ref{domin} below). If the sequence
$\{r_n(z_0)\}_{n\in \bN}$ diverges then  $\text{{\em supp}}\,\nu=\Xi_{f}$.
\item[(iii)] One has
$$\nu=\frac 1 \pi \frac {\partial y_{dom}}{\partial \bar z}\;
\Longleftrightarrow \; y_{dom}(z)=\int_{\bC}\frac {d\nu(\xi)}{z-\xi}.$$
\end{enumerate}
\end{theorem}

\begin{figure}[!htb]\label{fig1}
\centerline{\hbox{\epsfysize=5.0cm\epsfbox{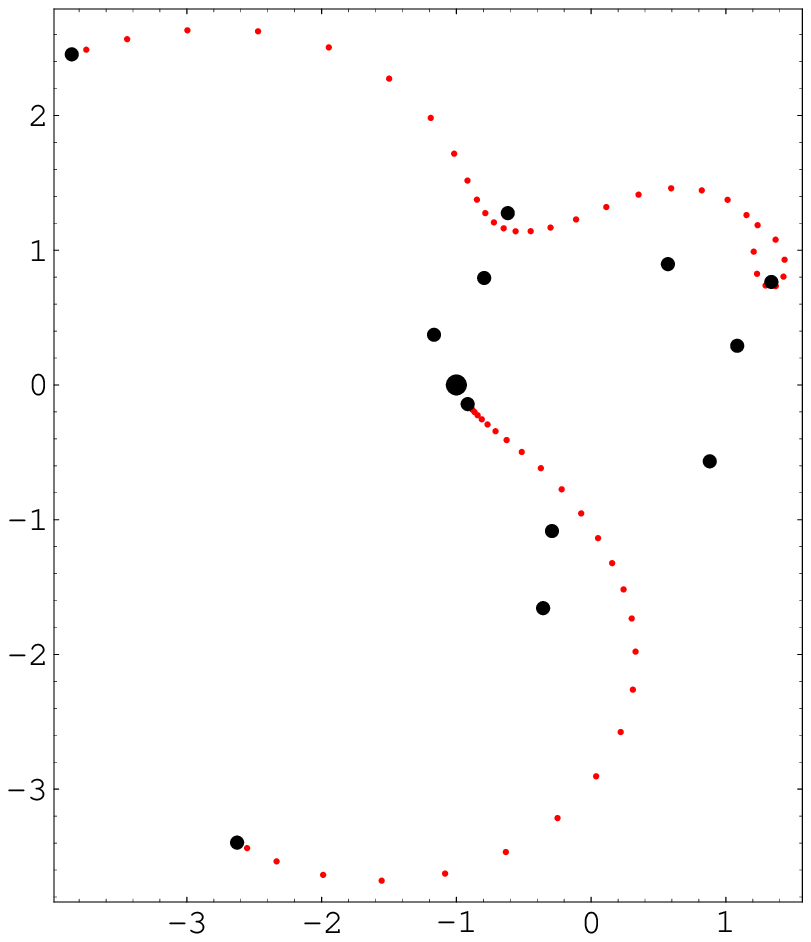}}
\hspace{0.8cm}\hbox{\epsfysize=5.0cm\epsfbox{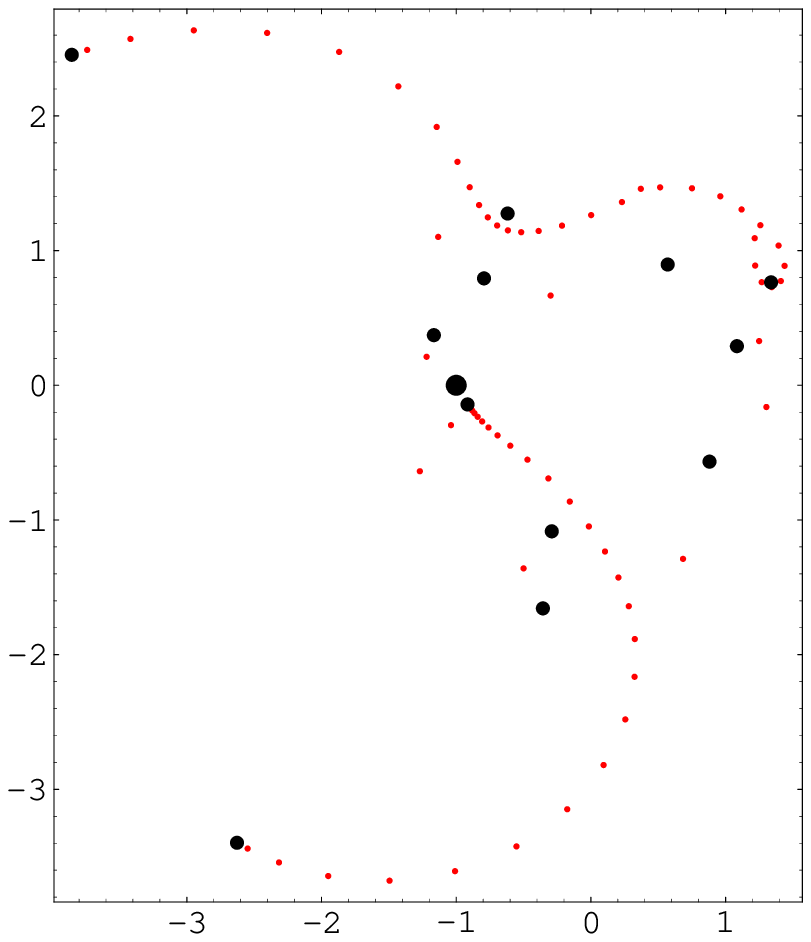}}}
\caption{Poles of $r_{31}(z)$ approximating the branch with maximal absolute
value of the algebraic function $f(z)$ with defining equation
$(z+1)y^3=(z^2+1)y^2+(z-5I)y+(z^3-1-I)$ for two different choices of
initial triples.}
\end{figure}

As a consequence of Theorems~\ref{th:Exist}--\ref{th:2} we obtain the following
asymptotic classification of the set of poles of the rational approximants
$\{r_n(z)\}_{n\in \bN}$ (see Figure~\ref{fig1}):

\begin{corollary}\label{cor1}
For any algebraic function $f(z)$ and any initial
        $k$-tuple $IN$ of rational functions the set of all poles of the family
$\{r_{n}(z)\}_{n\in \bN}$ splits
asymptotically into the following three types:
\begin{enumerate}
\item[(i)] The fixed part consisting of a subset of
$\Up_{f,IN}\setminus (\Xi_{f}\cup \Si_{f,IN})$;
\item[(ii)] The regular part tending asymptotically to the finite union
        of curves $\Xi_{f}$, that is, the induced equimodular discriminant
of~\eqref{eq:df1} (cf.~Theorem~\ref{th:Exist});
\item[(iii)] The spurious part tending to the (finite) set of slow growth
$\Si_{f,IN}$.
\end{enumerate}
\end{corollary}

The usage of the name ``spurious'' in Corollary~\ref{cor1} (iii) follows the
conventional way of describing poles of rational approximants to $f$
that do not correspond to
analytic properties of the function $f$. This terminology has been widely
adopted in
Pad\'e approximation, see e.g.~\cite{St5}. The rigorous definition of
spurious is rather technical, which is why we postpone it until \S 3 below.

\begin{expl1}
The two pictures show the poles of $r_{31}(z)$ for the initial triples
$p_{-2}(z)=p_{-1}(z)=0$, $p_{0}(z)=1$ and
$p_{-2}(z)=z^5+Iz^2-5$, $p_{-1}(z)=z^3-z+I$, $p_{0}(z)=1$, respectively (the
second initial triple was picked randomly).
The large fat point in these pictures is the unique pole $z=-1$ of $f$ which
in this case
coincides with $\Up_{f,IN}$
     since the chosen initial triples are polynomial.
The smaller fat points on both  pictures are the branching points of
$f(z)$.  One can show that the curve segments belonging to $\Xi_{f}$
can only end at these branching points, see \cite {Bi} and \cite
{BSh}. Note that the poles form the same pattern on both pictures. This
pattern is close to  the curve segments
constituting $\Xi_{f}$ and contains an additional
     nine isolated points in the right picture. These nine isolated points form
the set of slow growth
in this case; one can also
     check numerically that the latter set is indeed a subset of the
     finite discriminantal set $\cS$ described in the proof of
Lemma~\ref{descS}.
\end{expl1}

The final result of the paper settles a natural  analog of Nuttall's
Conjecture (see \S 3 and \cite {St5}) for our approximation scheme:

\begin{theorem}\label{th:Nut}
For any algebraic function $f(z)$ and any initial
        $k$-tuple $IN$ of rational functions there exists a finite upper bound
for the total number of spurious poles of the approximants
$\{r_n(z)\}_{n\in \bN}$ associated with the data $(f,IN)$.
\end{theorem}

As it was pointed out in~\cite{St5}, there is strong evidence that there
always exists some infinite subsequence of the sequence of diagonal Pad\'e
approximants to a given analytic function for which there are no spurious
poles in
a certain convergence domain. It is interesting to note that for the standard
choice of initial $k$-tuple $IN=\{0,0,\ldots,0,1\}$ in our approximation
scheme none of the rational functions in the resulting sequence of
approximants $\{r_n(z)\}_{n\in \bN}$  has spurious poles, see
Corollary~\ref{nospur}.

Let us conclude this introduction with a few remarks 
on the approximation scheme proposed above and some 
related topics. In the present setting the pointwise convergence of the 
sequence of rational approximants $\{r_n(z)\}_{n\in \bN}$ follows almost 
immediately from Theorem~\ref{th:ordrec} (cf.~Lemma~\ref{lm:basrec} below).
This is simply because our main recursion formula~\eqref{eq:rec} may actually 
be viewed as a recurrence relation with constant coefficients depending 
algebraically on a complex variable $z$. However, in many important 
cases one has to invoke quite different arguments in order to derive 
similar (pointwise) convergence results. For instance, the major challenge in 
the theory of finite difference equations is related to recurrence relations 
with varying coefficients. In this case the fundamental theorems of 
H.~Poincar\'e and O.~Perron (see, e.g., \cite[pp. 287-298]{Gel}) serve as 
prototypes for many new developments. Some interesting 
generalizations of the Poincar\'e and Van Vleck theorems were recently
obtained in \cite{Bu2,Bu3}. Typical applications of this sort of 
results pertain to the ratio asymptotics of various kinds 
of orthogonal polynomials and may be found in the well-known paper 
\cite{Ra}. We refer to \cite{Ba1} and \cite{Ba2} for more recent results 
in this direction. In particular, in {\em loc.~cit.~}the authors study the 
asymptotics of the zeros of a family of polynomials satisfying a usual 3-term
recurrence relation with varying and stabilizing complex coefficients and
prove that these zeros concentrate on an interval in the complex plane.
We continue this line of research in a forthcoming publication \cite{BSh} 
where we study the asymptotics of general families of functions satisfying
functional recurrence relations with varying and stabilizing
coefficients and obtain appropriate extensions of Poincar\'e's theorem.

This paper is organized as follows. In \S 2 we give all the
proofs, \S 3 contains a short description of Pad\'e
approximants and their spurious poles while \S 4 is devoted to
proving the $3$-conjecture of Egecioglu, Redmond and Ryavec. Finally, in \S 5
we compare our approximation scheme with Pad\'e approximation and discuss a
number of related topics and open problems.

\begin{ack}
Communication with Alan Sokal has been of great help in relating our
results to those of \cite{BKW1,BKW2}. The authors are obliged to the
anonymous referees for their constructive comments and for 
pointing out several references. 
\end{ack}

\section{Proofs and further consequences}\label{s2}

In order to prove Theorem~\ref{th:Exist} and Theorem~\ref{th:rate} we need
several preliminary results. Let us first fix the following notation:

\begin{notation}\label{new-not2}
Let $\Om$ be a domain in $\bC$ and $T:\Om\ni z\mapsto T(z)=(t_{ij}(z))\in
M_k(\bC)$ be a $k\times k$ matrix-valued map. Denote by
$\chi_{_{T}}:\Om\ni z\mapsto \chi_{_{T}}(\ze,z)\in \bC[\ze]$ the
characteristic
polynomial map associated with $T$, that is, for each $z\in \Om$ the
characteristic polynomial of $T(z)$ is given by the polynomial
$\chi_{_{T}}(\ze,z)$ in the variable $\ze$ with coefficients which are
complex-valued functions of $z$. Let further $\De_T\subset \Om$ be the
discriminant surface associated with $T$, i.e., the set of all
$z\in \Om$ for which the (usual) discriminant of the polynomial
$\chi_{_{T}}(\ze,z)$ vanishes. Finally, let $\Xi_T$
denote the equimodular discriminant associated with $T$, that is, the induced
equimodular discriminant of the family of polynomials
$\{\chi_{_{T}}(\ze,z)\mid z\in \Om\}$ (cf.~Definition~\ref{def7}).
\end{notation}

The following lemma is a consequence of well-known results on the
diagonalization of analytic
matrices which are based essentially on the implicit function theorem
(see~\cite[p.~106]{Ka}).

\begin{lemma}\label{l-diag}
If the map $T$ is analytic in $\Om$ then the eigenvalues $\la_i(z)$,
$1\le i\le k$, of $T(z)$ are analytic functions in $\Om\setminus \De_T$.
Furthermore, there exists an analytic map
$A:\Om\setminus \De_T\ni z\mapsto A(z)\in GL_k(\bC)$ such that
$$T(z)=A(z)\,\text{{\rm diag}}(\la_1(z),\ldots,\la_k(z))A(z)^{-1},\quad
z\in \Om\setminus \De_T.$$
The map $A(z)$ is uniquely determined in $\Om\setminus \De_T$ up to a
nonvanishing complex-valued analytic function $\gamma(z)$. Moreover,
its column vectors $A(z)^{(1)},\ldots,A(z)^{(k)}$ are eigenvectors of $T(z)$
with eigenvalues $\la_1(z),\ldots,\la_k(z)$, respectively, while the row
vectors of its inverse $A(z)_{(1)}^{-1},\ldots,A(z)_{(k)}^{-1}$ are
eigenvectors of the transpose $T(z)^t$ with eigenvalues
$\la_1(z),\ldots,\la_k(z)$,
respectively.
\end{lemma}

\begin{lemma}\label{l-forms}
Under the assumptions of Lemma~\ref{l-diag} the following holds:
\begin{enumerate}
\item[(i)] For any $\bx=(x_{k-1},\ldots,x_0)^t\in \bC^k$ there exist
complex-valued functions $\al_{ij}(z,\bx)$, $1\le i,j\le k$, which are analytic
in $\Om\setminus \De_T$ such that
\begin{equation}\label{form1}
\begin{split}
&T(z)^{n}\bx=\left(\sum_{j=1}^{k}\al_{1j}(z,\bx)\la_j(z)^n,\ldots,
\sum_{j=1}^{k}\al_{kj}(z,\bx)\la_j(z)^n\right)^t\\
&\text{for }z\in \Om\setminus \De_T\text{ and }n\in \bN.
\end{split}
\end{equation}
\item[(ii)] For $z\in \Om\setminus \Xi_T$ let $\la(z)=\la_1(z)$
denote the dominant eigenvalue of $T(z)$. There
exist $\bC^k$-valued analytic functions $\bu(z)=(u_1(z),\ldots,u_k(z))^t$ and
$\bv(z)=(v_1(z),\ldots,v_k(z))^t$ defined on $\Om\setminus \De_T$
such that
\begin{equation}\label{form2}
\begin{split}
&\bv(z)^t\bu(z)=1\text{ for }z\in \Om\setminus \De_T\text{ and}\\
&\frac{T(z)^n}{\la(z)^n}\rightrightarrows \bu(z)\bv(z)^t\text{ in }
\Om\setminus (\De_T\cup \Xi_T)\text{ as }n\to\infty.
\end{split}
\end{equation}
Thus $\left\{\frac{T(z)^n}{\la(z)^n}\right\}_{n\in \bN}$ is a sequence of
analytic matrices which converges locally uniformly in
$\Om\setminus (\De_T\cup \Xi_T)$ to $\bu(z)\bv(z)^t$ in the sup-norm on
$M_k(\bC)$.
\item[(iii)] In the notations of (i) and (ii) one has
\begin{equation}\label{form3}
\al_{i1}(z,\bx)=u_i(z)\sum_{j=1}^{k}v_j(z)x_{k-j},\quad 1\le i\le k,
\end{equation}
for any $\bx=(x_{k-1},\ldots,x_0)^t\in \bC^k$.
\end{enumerate}
\end{lemma}

\begin{proof}
Part (i) is a direct consequence of Lemma~\ref{l-diag}. To prove (ii) recall
the notations of Lemma~\ref{l-diag} and set $\bu(z)=A(z)^{(1)}$ and
$\bv(z)=A(z)^{-1}_{(1)}$, so that $\bv(z)^t\bu(z)=1$ for all
$z\in \Om\setminus \De_T$. From Lemma~\ref{l-diag} one gets
\begin{equation}\label{form4}
\frac{T(z)^n}{\la(z)^n}=A(z)\,\text{{\rm diag}}\!
\left(1,\left(\frac{\la_2(z)}{\la(z)}\right)^n,\ldots,
\left(\frac{\la_k(z)}{\la(z)}\right)^n\right)A(z)^{-1}
\end{equation}
for all $z\in \Om\setminus (\De_T\cup \Xi_T)$ and $n\in \bN$, which
immediately proves~\eqref{form2} since by assumption $|\la_j(z)|<|\la(z)|$ if
$2\le j\le k$. Finally, \eqref{form3} follows from (i) and (ii) by elementary
computations.
\end{proof}

Recall now the polynomials $P_i(z)$, $0\le i\le k$, from~\eqref{eq:df}
and define the rational functions
\begin{equation}\label{R}
R_i(z)=-\frac{P_i(z)}{P_k(z)},\quad 0\le i\le k-1.
\end{equation}
Choosing any initial $k$-tuple of rational functions
$IN=\{q_{0}(z),\ldots,q_{k-1}(z)\}$ one can rewrite recursion~\eqref{eq:rec} as
\begin{equation}\label{new-rec}
q_n(z)=\sum_{i=1}^{k}R_{k-i}(z)q_{n-i}(z)\text{ for }n\ge k
\end{equation}
or, equivalently,
\begin{equation}\label{mat-rec}
\begin{split}
&\bq_n(z)=T(z)^{n-k+1}\bq_{k-1}(z)\text{ for }n\ge k,\text{ where}\\
&\bq_m(z)=(q_m(z),q_{m-1}(z),\ldots,q_{m-k+1}(z))^t\in \bC^k
\text{ for }m\ge k-1\text{ and}\\
&T(z)=\begin{pmatrix} R_{k-1}(z) & R_{k-2}(z) & R_{k-3}(z)&
\ldots & R_0(z) \\
1 & 0 &0 & \ldots & 0\\
0 & 1 & 0 &\ldots &0 \\
\vdots & \vdots & \vdots & \ddots &\vdots \\
     0 &0 &\ldots &1 &0 \end{pmatrix}.
\end{split}
\end{equation}
We now apply Lemmas~\ref{l-diag} and~\ref{l-forms} to this specific case. It is
easy to see that the characteristic polynomial of $T(z)$ is given by
\begin{equation}\label{charpoly}
\chi_{_{T}}(\ze,z)=\ze^k-\sum_{i=0}^{k-1}R_{k-i-1}(z)\ze^{k-i-1}.
\end{equation}
Moreover, in view of Theorem~\ref{th:Exist}, Notation~\ref{new-not2}
and Lemma~\ref{l-forms} (ii) one has
\begin{equation}\label{nonv}
y_{dom}(z)=\la(z)\text{ for }z\in\bC\setminus (\Up_{f,IN}\cup \Xi_f),
\text{ where }\Xi_f=\Xi_T.
\end{equation}

\begin{remark}\label{r-nonv}
Identity~\eqref{nonv} implies in particular that $y_{dom}(z)$ is a nonvanishing
analytic function in $\bC\setminus (\Up_{f,IN}\cup \Xi_f)$.
\end{remark}

\begin{remark}\label{fin-D}
The discriminant $\De_T$ associated with the map $T(z)$ given
by~\eqref{mat-rec} is a finite subset of $\bC$ since all the entries of
$T(z)$ are rational functions of $z$ (cf.~\eqref{R}).
\end{remark}

It turns out that in our specific case the vector-valued functions $\bu(z)$
and $\bv(z)$
defined in Lemma~\ref{l-forms} (ii) can be computed explicitly:

\begin{lemma}\label{uv}
For $z\in \bC\setminus (\Up_{f,IN}\cup \Xi_f)$ let $\bu(z)$ and $\bv(z)$ be
eigenvectors with eigenvalue $\la(z)$ of $T(z)$ and $T(z)^t$, respectively.
There exist complex-valued functions $\ga_1(z),\ga_2(z)$
which are analytic and nonvanishing in $\bC\setminus (\Up_{f,IN}\cup \Xi_f)$
such that
\begin{equation*}
\begin{split}
&\bu(z)=\ga_1(z) (u_1(z),\ldots,u_k(z))^t,
\text{where }u_i(z)=\la(z)^{k-i},\,1\le i\le k,\\
&\bv(z)=\ga_2(z) (v_1(z),\ldots,v_k(z))^t,
\text{where }v_i(z)=\sum_{j=i}^{k}R_{k-j}(z)\la(z)^{i-j},\,1\le i\le k.
\end{split}
\end{equation*}
In particular, the maps $\bu(z)$, $\bv(z)$ and $\bu(z)\bv(z)^t$ are all
analytic in $\bC\setminus (\Up_{f,IN}\cup \Xi_f)$. Moreover, the functions
$\ga_1(z),\ga_2(z)$ may be chosen so that $\bv(z)^t\bu(z)=1$ for all
$z\in \bC\setminus (\Up_{f,IN}\cup \Xi_f)$.
\end{lemma}

\begin{proof}
Note that $\bu(z)$ and $\bv(z)$ are analytic $\bC^k$-valued functions
in $\bC\setminus (\Up_{f,IN}\cup \Xi_f)$ since $\la(z)$ is a simple eigenvalue
of $T(z)$ (and $T(z)^t$) for all $z$ in this set
(cf.~the proof of Lemma~\ref{l-diag}).
The explicit forms of $\bu(z)$ and $\bv(z)$ are easily obtained by
elementary manipulations with the matrix $T(z)$ defined in~\eqref{mat-rec}.
Note also that since $\la(z)$ is a root of the characteristic polynomial
$\chi_{_{T}}(\ze,z)$ of $T(z)$ one has
\begin{equation*}
\bv(z)^t\bu(z)=\ga_1(z)\ga_2(z)\sum_{m=0}^{k-1}(k-m)R_{m}(z)\la(z)^{m}=
\ga_1(z)\ga_2(z)\la(z)
\frac{\partial}{\partial \ze}\chi_{_{T}}(\ze,z)\Big|_{\ze=\la(z)}.
\end{equation*}
Therefore, if $\bv(z)^t\bu(z)=0$ then $\la(z)$ must be a multiple root of the
characteristic polynomial $\chi_{_{T}}(\ze,z)$ of $T(z)$. However, this
contradicts the assumption $z\in \Xi_f$. Thus $\bv(z)^t\bu(z)\neq 0$
whenever $z\in \Xi_f$, which proves the last statement of the lemma.
\end{proof}

Henceforth we assume that the complex-valued functions $\ga_1(z),\ga_2(z)$
defined in Lemma~\ref{uv} are chosen so that $\bv(z)^t\bu(z)=1$ for all
$z\in \bC\setminus (\Up_{f,IN}\cup \Xi_f)$.

\begin{remark}\label{non-0-u}
By Remark~\ref{r-nonv} all coordinates of $\bu(z)$ are
nonvanishing (analytic) functions in $\bC\setminus (\Up_{f,IN}\cup \Xi_f)$.
\end{remark}


\begin{lemma}\label{uc-q}
With the above notations one has
\begin{equation}\label{uc-q1}
\frac{q_n(z)}{\la(z)^n}\rightrightarrows
\ga_1(z)\ga_2(z)u_1(z)\sum_{i=1}^{k}v_i(z)q_{k-i}(z)\text{ in }
\bC\setminus (\Up_{f,IN}\cup \Xi_f)
\end{equation}
as $n\rightarrow \infty$.
\end{lemma}

\begin{proof}
Note that in view of recursion~\eqref{new-rec} all the functions
$\frac{q_n(z)}{\la(z)^n}$, $n\in \bN$, are analytic in
$\bC\setminus (\Up_{f,IN}\cup \Xi_f)$ while
Lemma~\ref{uv} implies that the same is true for the function in the
right-hand side of~\eqref{uc-q1}. From Lemma~\ref{l-forms} (i)-(iii) it
follows that
\begin{equation}\label{uc-q2}
\frac{q_n(z)}{\la(z)^n}\rightrightarrows
\ga_1(z)\ga_2(z)u_1(z)\sum_{i=1}^{k}v_i(z)q_{k-i}(z)\text{ in }
\bC\setminus (\Up_{f,IN}\cup \Xi_f\cup \De_T)
\end{equation}
as $n\rightarrow \infty$. Since $\De_T$ is a finite subset of $\bC$
(cf.~Remark~\ref{fin-D}) Cauchy's
integral formula shows that the local uniform convergence
stated in~\eqref{uc-q2} must actually hold even in an $\eps$-neighborhood of
the set $\De_T\setminus (\Up_{f,IN}\cup \Xi_f)$, which proves~\eqref{uc-q1}.
\end{proof}

To further simplify the notations set
\begin{equation}\label{gn-g}
g_n(z)=\frac{q_n(z)}{\la(z)^n}, n\in \bN, \text{ and }
g(z)=\ga_1(z)\ga_2(z)u_1(z)\sum_{i=1}^{k}v_i(z)q_{k-i}(z).
\end{equation}
We can complement Lemma~\ref{uc-q} with the following result:

\begin{lemma}\label{exp-q}
Given $0<\eps\ll 1$ let $\cV_{\eps}=\bCP\setminus V_{\eps}$, where $V_{\eps}$
is the $\eps$-neighborhood of $\Up_{f,IN}\cup \Xi_f$ in the spherical metric.
The rate of convergence of $g_n(z)\rightrightarrows g(z)$ in $\cV_{\eps}$ is
exponential, that is, there exist constants $\mathfrak{C}_{\eps}$ and
$\rho_{\eps}\in (0,1)$ such that $|g_n(z)-g(z)|\le
\mathfrak{C}_{\eps}\rho_{\eps}^n$ for all $z\in \cV_{\eps}$.
\end{lemma}

\begin{proof}
By Lemma~\ref{l-forms} (i) and~\eqref{form4} we get an exponential rate of
convergence for $g_n(z)\rightrightarrows g(z)$ in $\bCP\setminus U_{\eps}$,
where
$U_{\eps}$ is the $\eps$-neighborhood of $\Up_{f,IN}\cup \Xi_f\cup \De_T$ in
the spherical metric. Since $g_n(z)$, $n\in \bN$, and $g(z)$ are all analytic
functions in $\bC\setminus (\Up_{f,IN}\cup \Xi_f)$ and $\De_T$ is a finite
subset of $\bC$ we get the desired conclusion again by Cauchy's integral
formula.
\end{proof}

The next lemma gives a complete description of the
set of slow growth $\Si_{f,IN}$ that appears in Theorem~\ref{th:Exist}
(cf.~Definitions~\ref{df:slow} and~\ref{Sigma}):

\begin{lemma}\label{descS}
Let $IN=\{q_0(z),\ldots,q_{k-1}(z)\}$ be a fixed initial $k$-tuple of
rational functions as in~\eqref{new-rec}. Then
\begin{equation}\label{impl-sg}
\Si_{f,IN}=\left\{z\in \bC\setminus (\Up_{f,IN}\cup \Xi_f)\Bigg|
\sum_{i=1}^{k}\sum_{j=i}^{k}\la(z)^{i-j}q_{k-i}(z)R_{k-j}(z)=0\right\},
\end{equation}
so that in particular $|\Si_{f,IN}|<\infty$.
\end{lemma}

\begin{proof}
The fact that $\Si_{f,IN}$ is given by~\eqref{impl-sg} follows readily from 
Lemmas~\ref{l-forms}--\ref{uc-q} and 
Remark~\ref{non-0-u}. In order to show that $|\Si_{f,IN}|<\infty$ one can
produce an explicit upper bound for $|\Si_{f,IN}|$ in the following way: set
$$\mu_{_{T,IN}}(\ze,z)
=\sum_{i=1}^{k}\sum_{j=i}^{k}q_{k-i}(z)R_{k-j}(z)\ze^{i-j}$$
and let $S(z)$ denote the resultant of the polynomials
$\chi_{_{T}}(\ze,z)$ and $\mu_{_{T,IN}}(\ze,z)$ in the variable $\ze$. Let
further $\cS=\{z\in \bC\mid S(z)=0\}$. Since $S(z)$ is a rational
function and $\Si_{f,IN}\subseteq \cS$ it follows that
$|\Si_{f,IN}|\le |\cS|<\infty$.
\end{proof}

Both Theorem~\ref{th:Exist} and Theorem~\ref{th:rate} are immediate
consequences of Lemmas~\ref{uc-q}--\ref{descS} above. The same ideas yield
also a proof of Theorem~\ref{th:Nut}. Actually, we can produce an
asymptotically optimal bound for the number of spurious poles of the rational
functions $r_n(z)$. To do this let us introduce the following notation.

\begin{notation}\label{nota-spu}
Fix $\eps>0$ such that 
$\Si_{f,IN}\cap (\Up_{f,IN}\cup \Xi_f)_{\eps}=\emptyset$, where 
$(\Up_{f,IN}\cup \Xi_f)_{\eps}$ denotes the $\eps$-neighborhood of 
$\Up_{f,IN}\cup \Xi_f$ in the spherical metric. 
Given a sequence of rational functions $\{q_n(z)\}_{n\in \bN}$ as 
in~\eqref{new-rec} denote by $Z_{\eps}^{sp}(q_n)$ the set of zeros of $q_n$ 
which are 
contained in $\bC\setminus (\Up_{f,IN}\cup \Xi_f)_{\eps}$. 
The {\em total cardinality} 
$||Z_{\eps}^{sp}(q_n)||$ of $Z_{\eps}^{sp}(q_n)$ is the sum of the 
multiplicities of each 
of its elements as zeros of $q_n$. The {\em total cardinality} of 
$\Si_{f,IN}$ is $||\Si_{f,IN}||=\sum_{z\in \Si_{f,IN}}m(z)$, where 
$m(z)$ is the multiplicity of $z\in \Si_{f,IN}$ as a zero of the equation 
in~\eqref{impl-sg}.
\end{notation}

\begin{lemma}\label{l-Nut}
Let $IN=\{q_0(z),\ldots,q_{k-1}(z)\}$ be a fixed initial $k$-tuple of
rational functions. In the above notation there exists $N\in \bN$ such that
$||Z_{\eps}^{sp}(q_n)||\le ||\Si_{f,IN}||$ for $n\ge N$. Thus
$$||Z_{\eps}^{sp}(q_n)||\le 
\max(||Z_{\eps}^{sp}(q_0)||,\ldots,||Z_{\eps}^{sp}(q_{N-1})||,||\Si_{f,IN}||)$$
for all $n\in \bN$.
\end{lemma}

\begin{proof}
Recall~\eqref{gn-g} and note that since $|\Si_{f,IN}|<\infty$ it is enough to
show that for any $z_0\in \Si_{f,IN}$ there exists $N_0\in \bN$ and $\eps_0>0$
such that $g_n(z)$ has (at most) $m(z_0)$ zeros 
in $\{z\in\bC\mid |z-z_0|<\eps_0\}$ for $n\ge N_0$.
But this is an immediate consequence of the uniform convergence established in 
Lemma~\ref{uc-q} and Hurwitz's Theorem.
\end{proof}

An interesting consequence of the above results is that by using our scheme
one can always find
rational approximants $\{r_n(z)\}_{n\in \bN}$ that have no spurious poles:

\begin{corollary}\label{nospur}
If $f$ is an algebraic nonrational function and $IN$ is chosen to be the
standard initial $k$-tuple
$\{0,0,\ldots,1\}$ then $\Si_{f,IN}=\emptyset$. In particular,
the rational approximants $\{r_n(z)\}_{n\in \bN}$ corresponding to the
standard initial $k$-tuple have no spurious poles.
\end{corollary}

\begin{proof}
Indeed, if $\Si_{f,IN}\neq\emptyset$ and $z\in \Si_{f,IN}$ then
by~\eqref{charpoly} and~\eqref{impl-sg} the
corresponding dominant root $y_{dom}(z)=\la(z)$ has to satisfy
$$\la(z)^k-\sum_{i=0}^{k-1}R_{k-i-1}(z)\la(z)^{k-i-1}=0\,\text{ and }\,
\sum_{j=1}^{k}R_{k-j}(z)\la(z)^{1-j}=0,$$
which gives $\la(z)=0$. However, the latter identity cannot be fulfilled by
the dominant root of~\eqref{eq:rec} if $f$ is an algebraic nonrational
function.
\end{proof}

\begin{remark}
An analog of  Corollary~\ref{nospur}
for 3-term recurrence relations satisfied by certain biorthogonal
polynomials was proved in \cite[\S 3]{BSh}.
\end{remark}

\begin{remark}
Corollary~\ref{nospur} does not hold for arbitrarily chosen initial $k$-tuples
consisting of constant functions, as one can easily see by considering
e.g.~the initial $k$-tuple $IN=\{1,0,\ldots,1\}$.
\end{remark}

For the proof of Theorem~\ref{th:2} we need some additional notation and
results.

\begin{notation}\label{not1}
      Let $\mathfrak{Rec}_{k}$ be the $k$-dimensional complex linear space
         consisting of all $(k+1)$-term recurrence relations with constant
coefficients of the
form~\eqref{eq:Basic}. We denote by $\mathfrak{IN}_{k}$ the $k$-dimensional
complex linear space of all initial $k$-tuples $(u_{0},...,u_{k-1})$.
\end{notation}

Recall the notion of recursion of nondominant type introduced in
Definition~\ref{def5}.

\begin{definition}\label{domin}
A nondominant recurrence relation in $\mathfrak{Rec}_{k}$ with
initial $k$-tuple
$IN\in \mathfrak{IN}_k$ is said to be of {\em subdominant type} if the
following conditions
are satisfied. Let $\tau_1,\ldots,\tau_s$, $s\le k$, denote all distinct
spectral numbers with maximal absolute value and assume that they have
multiplicities $m_1,\ldots,m_s$, respectively. Then there exists a unique
index $i_0\in\{1,\ldots,s\}$ such that $m_i<m_{i_0}$ for
$i\in \{1,\ldots,s\}\setminus \{i_{0}\}$ and the initial $k$-tuple $IN$ is
{\em fast growing} in the sense that the degree of the polynomial $p_{i_0}$
in~\eqref{eq:leadasymp} corresponding to $\tau_{i_0}$ is precisely
$m_{i_0}-1$. The number $\tau_{i_0}$ is called the {\em dominant
spectral number} of this recurrence relation.
\end{definition}

The following lemma is a simple consequence of Theorem~\ref{th:ordrec}.

\begin{lemma}\label{lm:basrec}
In the above notation the following is true:
\begin{enumerate}
\item[(i)] The set of all slowly growing initial $k$-tuples
with respect to a given dominant recurrence relation in
$\mathfrak{Rec}_{k}$ is a complex
hyperplane $\SG_k$ in $\mathfrak{IN}_{k}$. The set $\SG_k$ is called the
{\em hyperplane of slow growth}.
\item[(ii)] For any dominant recurrence relation in
$\mathfrak{Rec}_{k}$ and any
fast growing
initial $k$-tuple in $\mathfrak{IN}_{k}$ the limit
$\lim_{n\to\infty}\frac{u_{n+1}}{u_{n}}$ exists and
         coincides with the dominant spectral number $\tau_{max}$, that is,
         the (unique) root of the characteristic equation~\eqref{eq:Char} with
         maximal absolute value.
\item[(iii)] Given a nondominant recurrence relation of
subdominant type in $\mathfrak{Rec}_{k}$ and a fast growing
initial $k$-tuple in $\mathfrak{IN}_{k}$ the limit
$\lim_{n\to\infty}\frac{u_{n+1}}{u_{n}}$ exists and
coincides with the dominant spectral number.
\item[(iv)] For any nondominant recurrence relation in $\mathfrak{Rec}_{k}$
which is not of subdominant type the set of initial
$k$-tuples for which
$\lim_{n\to\infty}\frac{u_{n+1}}{u_{n}}$
exists is a union of complex subspaces of $\mathfrak{IN}_{k}$ of
positive codimensions.
This union is called the {\em exceptional variety}.
\end{enumerate}
\end{lemma}

\begin{proof}
   In order to prove (i) notice that the coefficient $\kappa_{max}$ in
   Definition~\ref{df:slow} is a nontrivial
linear combination of the
   entries of the initial $k$-tuple $\{u_{0},\ldots,u_{k-1}\}$ with
coefficients depending
   on $\al_{1},\ldots,\al_{k}$. Therefore, the condition
   $\kappa_{max}=0$ determines a complex hyperplane $\SG_k$ in
   $IN_{k}$. One can easily see that the hyperplane of slow growth is the
direct sum of all Jordan blocks
   corresponding to the spectral numbers of a given
recurrence~\eqref{eq:Basic} other than the leading one.

The assumptions of part
   (ii) together with~\eqref{eq:leadasymp} yield
   $u_{n}=\kappa_{max}\tau_{max}^{n}+\ldots$ for $n\in \bZ_+$, where the dots
stand for the
   remaining terms in (\ref{eq:leadasymp}) corresponding to the
   spectral numbers whose absolute values are strictly smaller than
$|\tau_{max}|$.
Therefore, the quotient $\frac{u_{n+1}}{u_{n}}$ has a limit as $n\to\infty$
and this limit coincides with $\tau_{max}$, as required. By
   definition $\tau_{max}$ is a root of~\eqref{eq:Char}, which completes
   the proof of (ii). This last step can alternatively be carried out by
    dividing both sides of~\eqref{eq:Basic} by $u_{n-k}$ and then letting
$n\to \infty$. In view of Definition~\ref{domin} the same arguments show that
the assertion in (iii) is true as well.

For the proof of (iv) we proceed as follows. Take any
nondominant recurrence relation of the form~\eqref{eq:Basic} and let
$\tau_{1},\ldots,\tau_{r}$, $r\le k$, be all its distinct spectral numbers
with maximal absolute value. Thus $|\tau_i|=|\tau_{max}|$ if and only
$1\le i\le r$. Choose an initial
    $k$-tuple $IT=\{u_{0},\ldots,u_{k-1}\}$ and denote by $p_{1},...,p_{r}$ the
polynomials in~\eqref{eq:leadasymp} corresponding to
$\tau_{1},\ldots,\tau_{r}$, respectively,
    for the sequence $\{u_{n}\mid n\in \bZ_+\}$ constructed using the given
recurrence with initial $k$-tuple $IT$ as above. Assuming, as we may, that
our recurrence relation is nontrivial we get from~\eqref{eq:leadasymp} that
$\tau_i\neq0$ if $1\le i\le r$. We may further assume that the degrees
$d_1,\ldots,d_r$ of the polynomials $p_{1},...,p_{r}$, respectively,
satisfy $d_1\ge\ldots\ge d_r$.
A direct check analogous to the proof of part (ii)
   shows that if only the polynomial $p_{1}$ is nonvanishing then
$\lim_{n\to \infty}\frac{u_{n+1}}{u_{n}}=\tau_{1}$. If $r\ge 2$ and
$s\in \{2,\ldots,r\}$ is such that $p_{1},\ldots,p_{s}$ are all nonvanishing
polynomials among
$p_{1},\ldots,p_{r}$ then using again~\eqref{eq:leadasymp} we get
    $$\frac {u_{n+1}}{u_{n}}=\frac {p_{1}(n+1)+p_{2}(n+1)\left(\dfrac
    {\tau_{2}}{\tau_{1}}\right)^{n+1}+\ldots +p_{s}(n+1)\left(\dfrac
    {\tau_{s}}{\tau_{1}}\right)^{n+1}+o(1)}{p_{1}(n)+p_{2}(n)\left(\dfrac
    {\tau_{2}}{\tau_{1}}\right)^n+\ldots +p_{s}(n)\left(\dfrac
    {\tau_{s}}{\tau_{1}}\right)^n+o(1)}.$$
Since $\big|\frac{\tau_{i}}{\tau_{1}}\big|=1$ and $\tau_i\neq \tau_1$,
$2\le i\le s$, it follows that if $d_1=d_2$ then the expression in the
right-hand side has no limit as $n\to\infty$. Therefore, if such a limit
exists then $d_1>d_2$, which gives us a complex subspace of $IN_{k}$ of
(positive) codimension equal to $d_1-d_2$. Thus the exceptional variety is a
union of complex subspaces of $IN_{k}$ of (in general) different
codimensions.
\end{proof}

\begin{proof}[Proof of Theorem~\ref{th:2}]
To prove the first statement notice that under the assumptions
of the theorem the dominant root of the
asymptotic symbol equation~\eqref{eq:df1} is a  well-defined analytic
function in a sufficiently small neighborhood of $\Sigma_{f,IN}$.
Therefore, the residue distribution~\eqref{resid} associated to this dominant
root vanishes in a neighborhood of $\Sigma_{f,IN}$.
Thus the set $\Sigma_{f,IN}$ of slowly growing initial conditions can be
deleted from the support of the asymptotic ratio distribution $\nu$.

In order to show that under the nondegeneracy conditions stated in (ii)
the support of $\nu$ must coincide with $\Xi_{f}$ it is enough to prove that
the sequence $\{r_{n}(z)\}_{n\in \bN}$ diverges almost
everywhere in $\Xi_{f}$. This is actually an immediate consequence of the
fact that the nondegeneracy assumptions in (ii) imply that at almost every
point $z\in \Xi_{f}$ recursion~\eqref{eq:rec} is neither of dominant nor of
subdominant type
(cf.~Definitions~\ref{def5} and~\ref{domin}). Indeed, since the set
$\Si_{f,IN}$ of slow
growth is finite it follows from  Lemma~\ref{lm:basrec} (iv) that the
sequence $\{r_{n}(z)\}_{n\in \bN}$ diverges for almost all
$z\in \Xi_{f}$.

Finally, in order to prove (iii) notice that by
Definitions~\ref{def1}--\ref{c-trans}
any rational function $r$ is the Cauchy transform of its residue distribution
$\nu_{r}$ and that by Theorem~\ref{th:Exist} one has
$r_{n}(z)=\frac{q_{n}(z)}{q_{n-1}(z)}
\rightrightarrows y_{dom}(z)$ in $\bC\setminus \D_f$. Therefore, if
$\nu_{n}$ denotes the residue distribution of $r_n$ then the Cauchy transform
$\C_{\nu}$ of $\nu:=\lim_{n\to\infty}\nu_{n}$ equals $y_{dom}$. It then
follows from the standard relation
between a distribution and its Cauchy transform (see \eqref{res-dis}) that
$$\nu=\frac 1 \pi \frac {\partial \C_{\nu}}{\partial \bar z}=
\frac 1 \pi \frac {\partial y_{dom}}{\partial \bar z},$$
which completes the proof of the theorem.
\end{proof}

\begin{proof}[Proof of Corollary~\ref{cor1}]
The assertion is a direct consequence of the previous results since
Theorem~\ref{th:rate} implies that for any sufficiently small $\epsilon >0$
	    the complement $\Theta_{\epsilon}$ of the
	    $\epsilon$-neighborhood $\mathcal O_{\epsilon}$  of
	    $\D_{f}=\Xi_{f}\cup \Up_{f,IN}\cup \Si_{f,IN}$ in the
	    spherical metric contains
	    no poles of $r_{n}(z)$ for sufficiently large $n$.
Note in particular that by Theorem~\ref{th:2} (iii) the regular part of the
poles of the family $\{r_{n}(z)\}_{n\in \bN}$ must tend to $\Xi_{f}$ in such
a way that $\lim
	   \nu_{n}=\nu=\frac 1 \pi \frac {y_{dom}(z)}{\partial \bar z}$.
\end{proof}

\section{A brief excursion around Pad\'e approximation}\label{s3}

     In order to make this paper self-contained and to be able to
     compare our approximation scheme  with the classical
     Pad\'e scheme we give in this section a brief account of Pad\'e
      approximants. For a detailed overview of this subject we refer to
e.g.~\cite {BG-M} and \cite{St5}.
As in \S 1 we denote by $\cP_{k}$ the linear space of all
      complex polynomials of degree at most $k$.

\begin{definition}\label{def-p1}
Let $f(z)$ be a function analytic at infinity with power series expansion
$f(z)=\sum_{i=0}^{\infty}c_{i}z^{-i}$. For each pair $(m,n)\in \bN^2$ there
exist
      two polynomials $p_{m,n}\in \cP_{m}$ and $q_{m,n}\in
      \cP_{n}\setminus \{0\}$ such that
      \begin{equation} \label{eq:Pd}
      q_{m,n}\left(\frac 1 z\right)f(z)-p_{m,n}\left(\frac 1
z\right)=\mathcal O(z^{-m-n-1})
      \end{equation}
      as $z\to \infty$. The rational function
	  \begin{equation}\label{eq:Pade}
	  [m/n](z):=\frac {p_{m,n}(1/z)}{q_{m,n}(1/z)}
	  \end{equation}
	  is called the $(m,n)$-Pad\'e approximant to the function $f(z)$
	  (expanded at $\infty$).
	  \end{definition}

The notation introduced in Definition~\ref{def-p1} will be used throughout
this section. Pad\'e approximants can be seen as a generalization of
Taylor polynomials to the field of rational functions. Note that in view
of~\eqref{eq:Pd} these approximants are well-defined and exist uniquely for
any function $f(z)$ that has a power
series expansion at $\infty$. Pad\'e approximants are also closely related to
(generalized)
        orthogonal polynomials. As we shall now explain, this connection is
especially simple in
        the case of {\em diagonal Pad\'e approximants}, that is, approximants
of the
form
   $[n/n](z)$. To this end let us define the reverse
denominator polynomial $Q_{n}$ of $q_{n,n}$ by
	$$Q_{n}(z)=z^nq_{n,n}\left(\frac 1 z\right).$$

The following results may be found in \cite[Lemma 2.2]{St5}.

\begin{lemma}
	    A polynomial $q_{n,n}\in
      \cP_{n}\setminus \{0\}$ is the denominator of the Pad\'e
	    approximant $[n/n]$ to $f$ if, and only if, the reverse polynomial
	    $Q_{n}$ of $q_{n,n}$ satisfies the (generalized) orthogonality
relation
	    \begin{equation}\label{eq:orth}
	    \int_{C}\zeta^iQ_{n}(\zeta)f(\zeta)d\zeta=0,\quad
i\in\{0,\ldots,n-1\},
	    \end{equation}
      where $C$ is a closed contour containing all the
      singularities of $f(z)$ (so that $f(z)$ is analytic on $C$ and in its
      exterior).
	    \end{lemma}

Condition~\eqref{eq:orth} is a (generalized)
		orthogonality relation for the polynomial family
$\{Q_{n}(z)\}_{n\in \bN}$
		on the contour $C$ with respect to
		the complex-valued measure $f(\zeta)d\zeta$.
		The fact that the quadratic form
		$$(g,h)_{f}:=\int_{C}g(\zeta)h(\zeta)f(\zeta)d\zeta$$
		diagonalized by the family $\{Q_{n}(z)\}_{n\in \bN}$ is
in general neither Hermitian nor
		positive-definite accounts for the rather
chaotic behavior of
		the zero loci of the polynomials $Q_{n}(z)$, $n\in \bN$, see
\cite{St2,St3,St5} and the results below.

          \begin{remark}
Note that the generalized orthogonality relation~\eqref{eq:orth} and Favard's
theorem imply that if the polynomials in the family
          $\{Q_{n}(z)\}_{n\in \bN}$ are assumed to be monic then they
satisfy a three-term recurrence
          relation of the form
          \begin{equation} \label{eq:3term}
	  zQ_{n}(z)=Q_{n+1}(z)+\al_{n}Q_{n}(z)+\be_{n}Q_{n-1}(z)
          \end{equation}
	  with complex-valued coefficients $\al_{n}$ and $\be_{n}$, see
\cite {Ch}.
	  Apparently, the chaotic behavior of the zero loci of the polynomials
	  $Q_{n}(z)$, $n\in \bN$, is inherited by the families of
	  coefficients $\{\al_{n}\}_{n\in \bN}$ and $\{\be_{n}\}_{n\in \bN}$,
as these do not necessarily converge, see \S 5.2.3.
\end{remark}

Probably the most important general theorem that applies to functions
meromorphic in the complex plane is that of Nuttall-Pommerenke asserting
that if $f$ is meromorphic in $\bC$ and analytic at $\infty$ then the
sequence of diagonal
	Pad\'e approximants $\{[n/n](z)\}_{n\in \bN}$ converges to $f(z)$ in
planar measure (cf.~\cite{BG-M,St1}).
The convergence of the diagonal
	Pad\'e approximants to $f(z)$ in the case when $f(z)$ is a
	(multi-valued) algebraic function was studied in
	detail in \cite {St1}. As pointed out in {\em loc.~cit.}, the following
important result was
	proved already in \cite {Du}.

	\begin{theorem} \label{th:du}
	    Let $f(z)$ be a locally meromorphic function in some domain
	    $\bCP^1\setminus E$, where $E$ is a compact subset of $\bC$
with $\text{{\em cap}}(E)=0$.
	    There exists a domain $D_{f}\subseteq \bCP^1$ such that the
sequence of  diagonal Pad\'e approximants $\{[n/n](z)\}_{n\in \bN}$ converges
   in capacity to
	    $f(z)$ in $D_{f}$ but it does
	    not converge in capacity to $f(z)$ in any domain $\widetilde
	    D\subseteq \bCP^1$ with
$\text{{\em cap}}(\widetilde D\setminus D)>0$.
\end{theorem}

The domain $D_{f}$ in Theorem~\ref{th:du} is uniquely
determined up to a set of zero
capacity and is called
	    the {\em convergence domain} of $f(z)$.  (Concerning the
notion of (logarithmic) capacity
         and convergence in capacity one may consult e.g.~\cite{ST}.)
	 The following improvement of Theorem~\ref{th:du} was obtained in
	 \cite {St5} by building on the results of \cite {St6}.

	 \begin{theorem} \label{th:qdiff}
Let $f(z)$ be a locally meromorphic function in some domain
	    $\bCP^1\setminus E$, where $E$ is a compact subset of $\bC$
with $\text{{\em cap}}(E)=0$. The convergence domain $D_{f}$ of $f(z)$ is
uniquely determined (up to a set of
	     capacity zero) by the following two conditions:
\begin{enumerate}
\item[(i)] $D_{f}$ is a subdomain of $\bCP^1$ containing $\infty$ such
	     that $f(z)$ has a single-valued meromorphic continuation
	     throughout $D_{f}$;
\item[(ii)] $\text{{\em cap}}([\bCP^1\setminus
	     D_{f}]^{-1})=\inf_{\widetilde
	     D}\text{{\em cap}}([\bCP^1\setminus \widetilde D]^{-1})$, where
the
	     infimum ranges over all domains $\widetilde D\subseteq \bCP^1$
	     satisfying (i). (Here for a given set $S\subset \bCP^1$ we denote
	     by $S^{-1}$ the set of all $z\in \bCP^1$ such that $z^{-1}\in
	     S$.)
\end{enumerate}
                 \end{theorem}

		 In other words, Theorem~\ref{th:qdiff} asserts that
the convergence domain $D_{f}$ of $f(z)$ is  characterized
by the property that its boundary $\partial D_{f}$ has minimal capacity.

	\begin{definition}\label{def-p2}
Given a function $f(z)$ analytic near $\infty$ we say that a compact
	set $K\subseteq \bCP^1$ is {\em admissible} for $f(z)$ if $f(z)$ has a
	single-valued continuation throughout $\bCP^1\setminus K$.
	  \end{definition}

Using the terminology introduced in Definition~\ref{def-p2} one may rephrase
Theorem \ref{th:qdiff} as follows: given a locally meromorphic function
$f(z)$ in some domain
	    $\bCP^1\setminus E$, where $E$ is a compact subset of $\bC$
with $\text{cap}(E)=0$, there exists
an unique admissible compact set
$K_{0}\subseteq \bCP^1$ such that
$$\text{cap}(K_{0})=\inf_{K}\text{cap}(K),$$
where the infimum is taken over all admissible compact sets
$K\subseteq \bCP^1$.

Let $g_{D}(z,w)$ denote the Green's function of a
	      domain $D\subset \bCP^1$ with $\text{cap}(\partial D)>0$.
The following result may be found in \cite[Part III,
		  Theorem 1 and Lemma 5]{St6}.

		  \begin{theorem}
In the above notation
the minimal set
		  $K_{0}$ is the union of a compact set $E_{0}$ with analytic
		  Jordan arcs. The set $E_{0}$ has the property that
the analytic
		  continuation of $f(z)$ in $\bCP^1\setminus K_{0}$ has a
		  singularity at every boundary point of $E_{0}$ while the
analytic arcs
		  in $K_{0}\setminus E_{0}$ are trajectories of the quadratic
		  differential $T_{f}(z)dz^2$ given by
		  $$T_{f}(z)=\left\{\frac {\partial}{\partial
z}g_{\bCP^1\setminus
		  K_{0}}(z,\infty)\right \}^2.$$
Moreover, if $f(z)$ is
		      algebraic then $T_{f}(z)$ is a
rational function.
\end{theorem}

		      \begin{theorem}\label{th:9}
Let $f(z)$ be an algebraic
nonrational function which is
			  analytic near $\infty$. If all
singularities of $f(z)$ are
			  contained in a set $W=\{w_{1},\ldots,
w_{n}\}\subset \bC$ then there
			  exist a subset $W'$ of $W$ with $|W'|=n'\ge 2$
and a set
			  $Z\subset \bC$ with (not necessarily distinct) $n'-2$
			  points such that
			  $$T_{f}(z)=\frac {\Pi_{z_{j}\in
Z}(z-z_{j})}{\Pi_{w_{j}\in
			  W'}(z-w_{j})}=c\left\{\frac
{\partial}{\partial z}g_{\bCP^1\setminus
		  K_{0}}(z,\infty)\right \}^2.$$
		  Here $c$ is an appropriate complex constant and $z\in
		  \bCP^1\setminus K_{0}$, where $K_{0}$ is the union of
trajectories of the quadratic
		  differential $T_{f}(z)dz^2$ and a set of isolated points
contained in $W$.
\end{theorem}

Theorem~\ref{th:9} shows that most of the poles of the diagonal
Pad\'e approximants $[n/n](z)$ tend
to $K_{0}$ when $n\to \infty$. However, not all of them do! Indeed, even
in the case when $f(z)$ is algebraic one cannot hope for a better convergence
type than convergence in capacity (which is only slightly stronger than
convergence almost everywhere with respect to Lebesgue measure). As shown by
the following example of H.~Stahl, uniform convergence in
$\bCP^1\setminus K_{0}$ fails in a rather dramatic way.

\begin{example}
Consider the function
$$f(z)=\frac{(z-\cos {\pi
\al_{1}})(z-\cos{\pi \al_{2}})}{\sqrt{z^2-1}}-z +(\cos {\pi
\al_{1}}+\cos {\pi \al_{2}}),$$
where $1,\al_{1},\al_{2}$ are
rationally independent numbers. Then the reverse denominators
$\{Q_{n}(z)\}_{n\in \bN}$
satisfy the generalized orthogonality relation
$$\int_{-1}^{1}Q_{n}(x)Q_{m}(x)\frac{(x-\cos {\pi\al_{1}})(x-\cos{\pi
\al_{2}})}{\pi\sqrt{1-x^2}} dx=0\text{ for }m\neq n.$$
Even though each polynomial $Q_n(z)$ has at most two zeros  outside
$[-1,1]$, the zeros of $Q_n(z)$ cluster everywhere in $\bC$
as $n\rightarrow \infty$ (see \cite[Example 2.4]{St5}).
        \end{example}

        A naive definition of {\em spurious poles} of Pad\'e approximants
        would be that these are the poles which are  not located near
        $K_{0}$. A more rigorous definition is
        as follows (cf.~\cite[Definition 4.1]{St5}).

\begin{definition}\label{def-spur}
Let $f$ be a function satisfying the conditions of
Definition~\ref{def-p1}. Let further $\N\subseteq \mathbb N$ be an infinite
sequence and $\{[n/n]\}_{n \in \N}$ be the corresponding subsequence
        of diagonal Pad\'e approximants to $f$. We define {\em
        spurious poles} in two different situations:
\begin{enumerate}
\item[(i)] Assume that for each $n\in \N$ the approximant $[n/n]$ has a
pole at $z_{n}\in \bCP^1$ such that $z_{n}\to z_{0}$ as $n\to \infty$,
$n \in \N$.
       If $f$ is analytic at $z_{0}$ and the
       approximants $[n/n]_{n \in \N}$ converge in capacity to $f$ in
       some neighborhood of $z_{0}$ then the poles of the approximants $[n/n]$
at $z_{n}$, $n\in \N$, are
       called {\em spurious}. In case $z_{0}=\infty$ the convergence
       $z_{n}\to z_{0}$ has to be understood  in the spherical metric.
\item[(ii)] Let the function $f(z)$ have a pole of order $k_{0}$ at
       $z_{0}$ and assume that for each $n\in \N$ the total order of poles of
the approximant $[n/n]$
       near $z_{0}$ is $k_{1}=k_{1,n}>k_{0}$. Assume further
       that $[n/n]$ has poles at $z_{n,j}$, $j\in \{1,\ldots,m_{n}\}$, of total
       order $k_{1,n}$ and that for any selection of
$j_{n}\in \{1,\ldots,m_{n}\}$ one has $z_{n,j_{n}}\to z_{0}$ as
$n\to \infty$, $n\in \N$ . Then poles
       of order $k_{1}-k_{0}$ out of all poles of the approximants
       $[n/n]$ near $z_{0}$, $n\in \N$, are called {\em spurious}.
\end{enumerate}
\end{definition}

With this definition we can finally formulate the
Baker-Gammel-Wills (Pad\'e) Conjecture and some related questions.

\begin{pade}
If the function $f$ is meromorphic
	    in the unit disk $\bD$ and analytic at $0$, and
if $\{[n/n]\}_{n\in \bN}$ denotes the
	    sequence of Pad\'e approximants to $f$ expanded at the origin,
	    then there exists an infinite subsequence $\N\subseteq \bN$
	    such that
	    $[n/n](z)\rightrightarrows f(z)$ in $\bD\setminus
\{\text{poles of f}\}\text{ as }n\to
\infty,\,n\in\N$.
\end{pade}

The Pad\'e Conjecture is widely regarded as the
main challenge in the theory of Pad\'e approximation. The
importance of this conjecture is due to the homographic invariance property
of diagonal Pad\'e approximants and its many consequences, see \cite{BG-M}.
A major breakthrough was recently made by D.~Lubinsky~\cite{Lu} who
showed that the general version of the Pad\'e Conjecture stated above is false.
Subsequently, V.~Buslaev~\cite{Bu} constructed counterexamples to this
conjecture for some special types of algebraic (hyperelliptic) functions.
However, the important question whether the
Pad\'e Conjecture could hold for certain classes of algebraic functions is
still open. A key step in this direction would be a deeper understanding of
the asymptotic distribution of spurious poles of diagonal Pad\'e approximants.
To this end J.~Nuttall proposed the following conjecture (cf.~\cite{St5}):

\begin{nutt}
Let $f$ be an algebraic
function which is analytic at $\infty$. Then there exists an upper
	bound for the number of spurious poles (in the sense of total
	order) for all Pad\'e approximants  $[n/n]$, $n\in \bN$.
\end{nutt}

Note that Theorems~\ref{th:Exist} and~\ref{th:Nut} show that
natural analogs
of both the Pad\'e Conjecture and Nuttall's
Conjecture are always true
for the rational approximants to algebraic functions constructed by using our
approximation scheme.

\section{Application: The $3$-conjecture of Egecioglu,
Redmond and Ryavec}
The following interesting conjecture appeared in~\cite{ERR} -- where it was
nicknamed the 3-conjecture -- in connection
with counting problems for various combinatorial objects such as alternating
sign matrices with vertical symmetries and $n$-tuples of nonintersecting
lattice paths.

\begin{conjecture}\label{3Con}
All polynomials belonging to the sequence $\{p_{n}(z)\}_{n\in \bN}$
defined by the $4$-term recursion
\begin{equation}\label{eq:3conj}
\begin{split}
& p_{n}(z)=zp_{n-1}(z)-Cp_{n-2}(z)-p_{n-3}(z),\text{ where}\\
& p_{-2}(z)=p_{-1}(z)=0,\,p_{0}(z)=1\text{ and }C\in \bR,
\end{split}
\end{equation}
have real zeros if and only if $C\ge 3$. Moreover, if $C>3$ then the zeros of
$p_{n+1}(z)$ and $p_{n}(z)$ are interlacing for all $n\in \bN$.
\end{conjecture}

Theorem 1 of~\cite{ERR} proves the sufficiency part as well as the last
statement in the above conjecture. We shall now fully settle the 3-conjecture
by showing that $C\ge 3$ is indeed a necessary condition for the validity of
this conjecture.

\begin{proposition}\label{pro-nec}
If $C<3$ and $\{p_{n}(z)\}_{n\in \bN}$ is the sequence of polynomials defined
by recursion~\eqref{eq:3conj} then there exists $N\in \bN$
such that not all zeros of $p_{N}(z)$ are real.
\end{proposition}

The proof of Proposition~\ref{pro-nec} is based on the ideas developed in
the previous sections and goes as follows. We study the behavior of the
        branching points of the curve $\Gamma_{C}$ given by solutions $t$ of
the asymptotic symbol equation
\begin{equation}\label{aseq}
     t^3-zt^2+Ct+1=0
\end{equation}
considered as a branched covering of the
        $z$-plane (cf.~Definition~\ref{def4}). More precisely, we show that
for any $C<3$ the
curve $\Gamma_{C}$ has two
        complex conjugate branching points and the induced equimodular
        discriminant $\Xi_{C}$ of~\eqref{aseq} is a real
semialgebraic curve ending at these
        branching points (cf.~\cite{Bi,BSh}). Therefore, if $\epsilon >0$ there
exists some
        sufficiently large $N\in \bN$ such that for all $n\ge N$ the polynomial
$p_{n}(z)$ has (complex conjugate) roots lying in the
$\epsilon$-neighborhoods of these branching  points.

\begin{figure}[!htb]
\centerline{\hbox{\epsfxsize=3.5cm\epsfbox{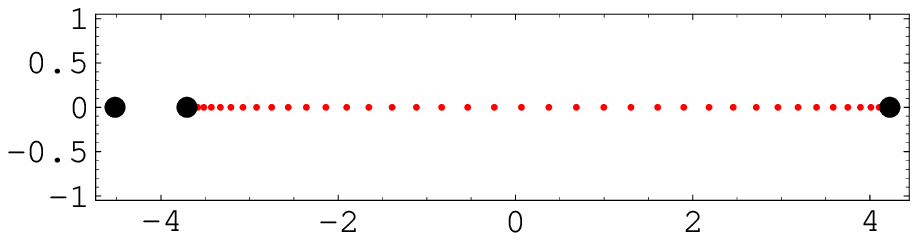}}
\hspace{0.5cm}\hbox{\epsfxsize=3.5cm\epsfbox{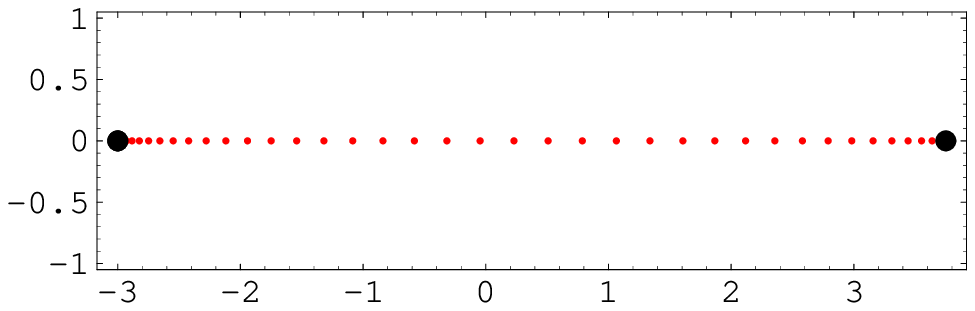}}
\hspace{0.5cm}\hbox{\epsfxsize=3.5cm\epsfbox{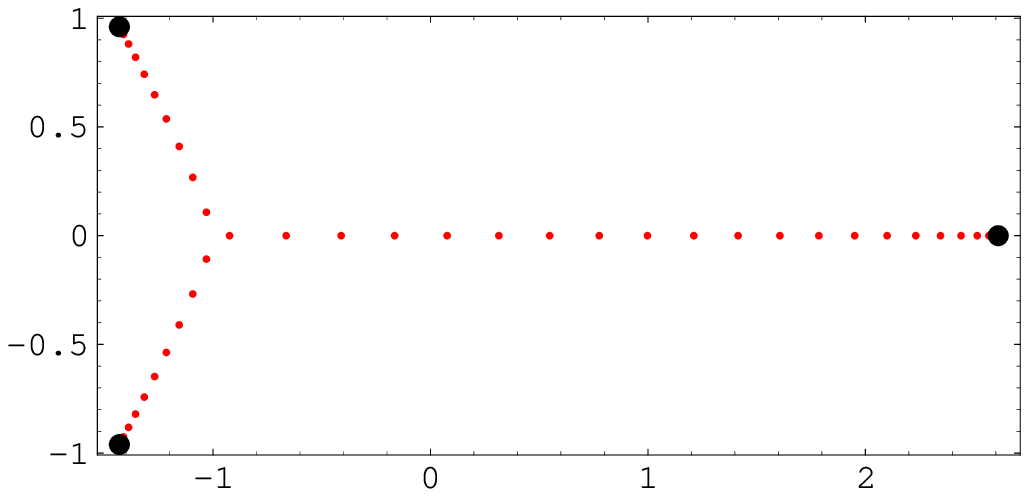}}}
\centerline{\hbox{\epsfysize=2.5cm\epsfbox{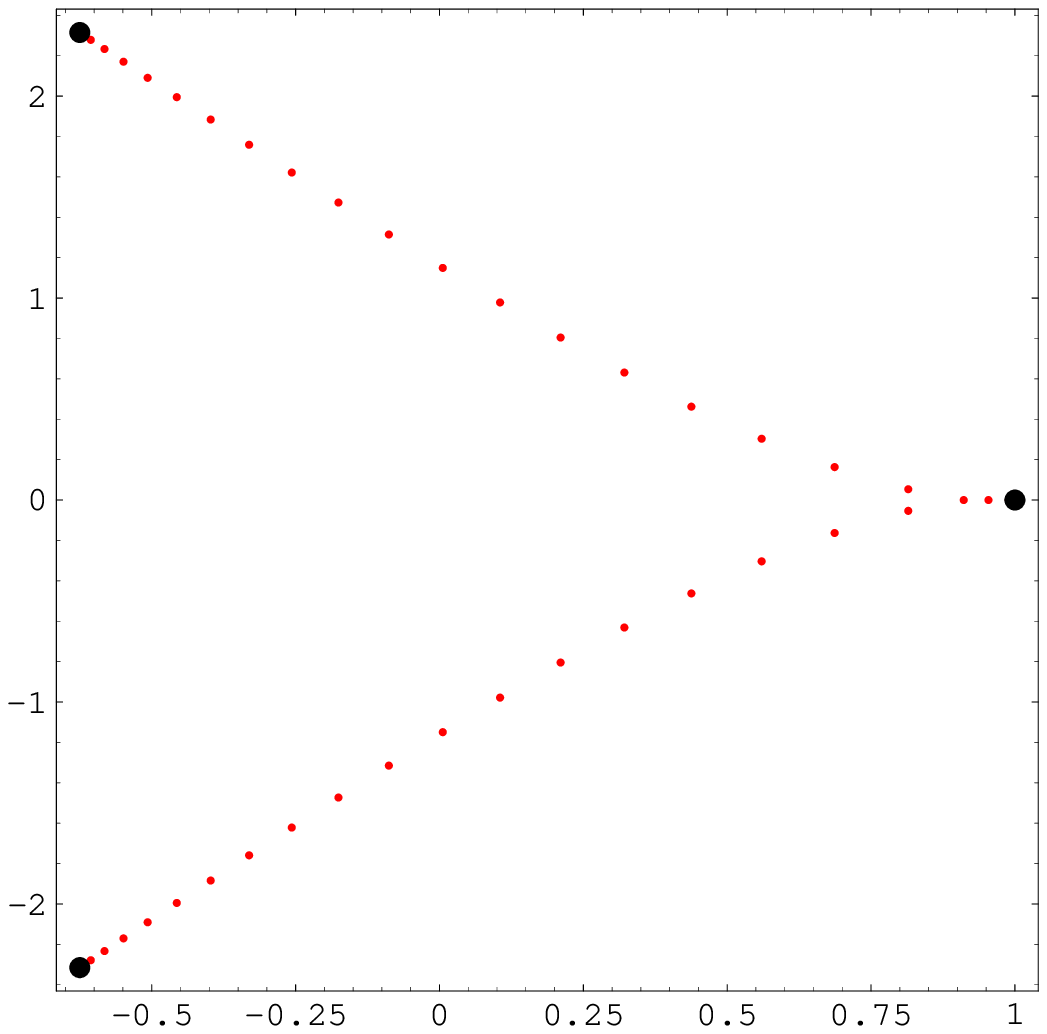}}
\hspace{0.5cm}\hbox{\epsfysize=2.5cm\epsfbox{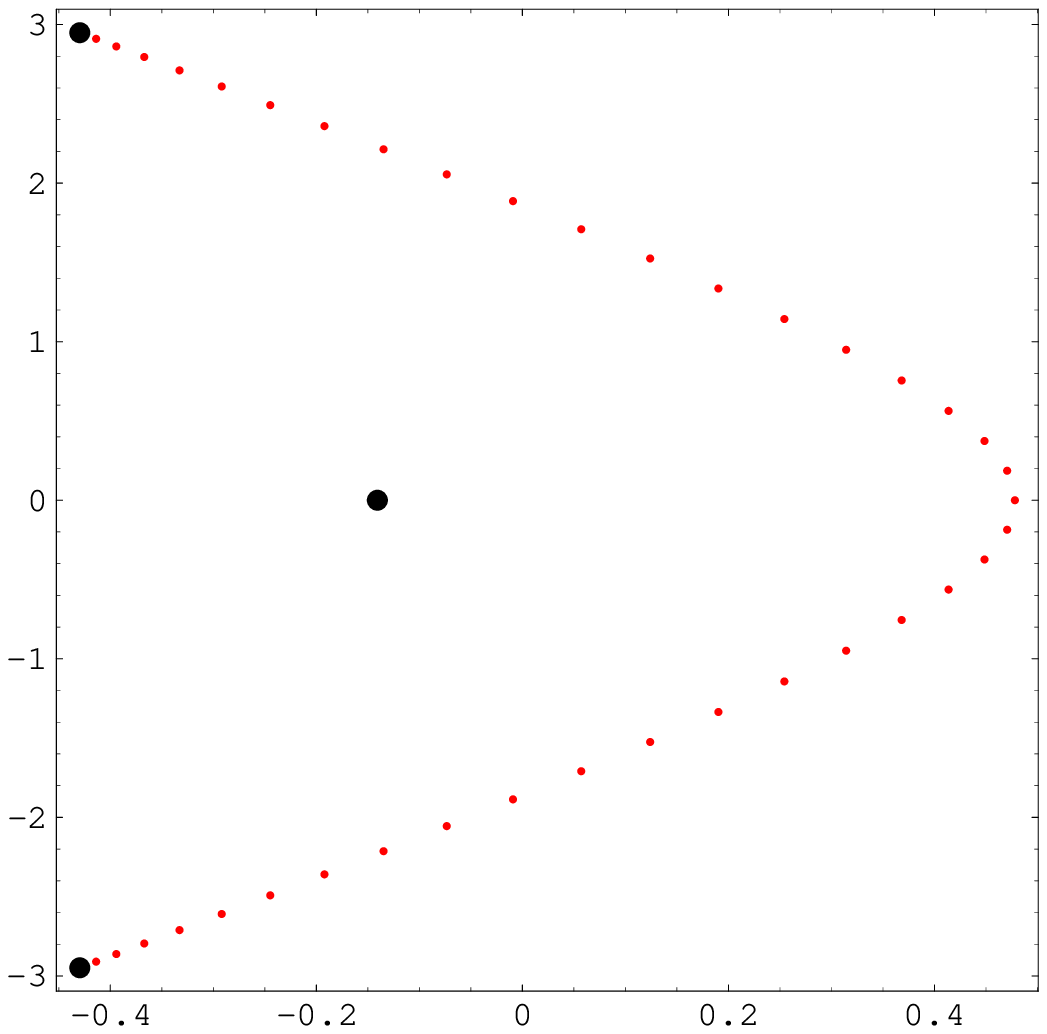}}}
\caption{Zeros of $p_{41}(z)$ satisfying~\eqref{eq:3conj} for
$C=4,\,3,\,1,\,-1,\,-2$.}\label{fig2}
\end{figure}

\begin{expl2}
The five pictures show the qualitative changes of the induced equimodular
discriminant $\Xi_C$ and the branching points of the curve $\Gamma_C$ when the
parameter $C$ runs from $\infty$ to $-\infty$. The fat points are the
branching points of $\Gamma_C$ while the thin points
illustrate the zeros of the polynomial $p_{41}(z)$. The first picture depicts
the typical behavior for $C>3$. In this case the branching points of
$\Gamma_C$ are all real and the equimodular discriminant $\Xi_C$ is the
segment joining the two rightmost branching points. If $C=3$ then the two
leftmost branching points of $\Gamma_C$ coalesce. For $-1<C<3$ there are one
real and two complex conjugate branching points of $\Gamma_C$ and the induced
equimodular discriminant $\Xi_C$ is the $Y$-shaped figure joining these three
points. If $C=-1$ then the triple point occurring in the previous case
coincides with the (unique) real branching point of $\Gamma_C$. Finally, if
$C<-1$ then one has a stable configuration consisting of the equimodular
discriminant $\Xi_C$ joining the two complex conjugate branching points of
$\Gamma_C$ and the remaining (real) branching point of $\Gamma_C$ lying
outside $\Xi_C$.
\end{expl2}

Straightforward computations using Mathematica$^{\textrm{TM}}$ yield the
following result:

\begin{lemma}\label{l-3conj1}
The branching points of $\Gamma_{C}$ satisfy the equation
\begin{equation}\label{eq:bran}
        4z^3+C^2z^2-18Cz-27-4C^3=0.
\end{equation}
The discriminant of the polynomial in $z$ occurring in the left-hand side
of~\eqref{eq:bran} is given by
\begin{equation}\label{eq:disc}
       64(C-3)^3(C^2+3C+9)^3.
\end{equation}
\end{lemma}

\begin{lemma}\label{l-3conj2}
The behavior of the induced equimodular
discriminant $\Xi_{C}$ near a branching point of $\Gamma_{C}$
changes if and only if the branching point becomes degenerate, i.e.,
the absolute value of the double root of the asymptotic symbol equation
corresponding to the branching point coincides with the absolute
value of the remaining root. This happens precisely for those values
of $C$ for which there exist $z,\,\tau,\,\theta$ such that the relation
\begin{equation}\label{eq:3-conj}
        t^3-zt^2+Ct+1=(t-\tau)^2(t-\tau e^{i\theta})
\end{equation}
is satisfied  for all $t\in \bC$. If $C$ is real this occurs only for $C=3$
and $C=-1$.
\end{lemma}

\begin{proof}
The first part of the lemma is obvious and can be found in e.g.~\cite {Bi}.
To prove the statement
about the values of $C$ consider the system for the coefficients
of~\eqref{eq:3-conj}:
\begin{equation*}
\begin{cases}
z\!\!\!\!&=\tau(2+e^{i\theta})\\
C\!\!\!\!&=\tau^2(1+2e^{i\theta})\,.\\
1\!\!\!\!&=-\tau^3e^{i\theta}
\end{cases}
\end{equation*}
The last equation gives $\tau=-e^{i(-\theta+2k\pi)/3}$, where $k\in \{0,1,2\}$.
Substituting this expression in the second equation one gets
$$C=e^{-2i(\theta+4k\pi)/3}+2e^{i(\theta+4k\pi)/3}.$$
Elementary computations now show that $C$ is real if and only if
$\theta=(3m-4k)\pi$ for some $m\in \bZ$, which gives either
$C=3$ (in which case $z=-3$ and $\tau=-1$) or $C=-1$ (in which case $z=1$ and
$\tau=1$).
\end{proof}

\begin{proof}[Proof of Proposition~\ref{pro-nec}]
Let us first show that it suffices to check that the following three
conditions hold for any value of $C\in (-\infty,3)$:
\begin{enumerate}
\item[(i)]  the curve $\Gamma_{C}$ has two (distinct) complex conjugate
branching points;
\item[(ii)] at each of these two complex branching points of $\Gamma_{C}$ the
absolute value of the double root of the corresponding asymptotic symbol
equation~\eqref{aseq} exceeds the absolute value of the remaining root;
\item[(iii)] at a generic point $z_0\in \Xi_{C}$ the roots
of the asymptotic symbol equation~\eqref{aseq}
with maximal absolute value have distinct arguments and the sequence
$\left\{\frac{p_n(z_0)}{p_{n-1}(z_0)}\right\}_{n\in \bN}$ diverges.
\end{enumerate}
Indeed, if (i) and (ii) are true then the induced equimodular discriminant
$\Xi_C$ must necessarily pass through (actually, end at) the two complex
conjugate branching points of $\Gamma_{C}$ (cf.~\cite{Bi,BSh}). In particular,
there exist some
nonisolated point $\ze\in \Xi_C$ and $\eps>0$ such that
$\{z\in \bC\mid |z-\ze|<\eps\}\subset \bC\setminus \bR$. Moreover,
if (iii) holds then Theorem~\ref{th:2} (ii) implies that
$\Xi_C=\text{supp}\,\nu$, where $\nu$ denotes the asymptotic ratio
distribution of the family $\{p_n(z)\}_{n\in \bN}$. Therefore, there exists
some sufficiently large $N\in \bN$ such that for all $n\ge N$ the polynomial
$p_{n}(z)$ has at least one root in
$\{z\in \bC\mid |z-\ze|<\eps\}\subset \bC\setminus \bR$, as required.

To verify conditions (i)-(iii) note first that if $C<3$ then~\eqref{eq:disc}
implies that all three branching points of $\Gamma_C$ are necessarily
distinct. Since~\eqref{eq:bran} is a polynomial equation with real
coefficients depending continuously on $C$ it is enough to check
that~\eqref{eq:bran} has (distinct) complex conjugate roots for an arbitrarily
chosen value of $C$ in the interval $(-\infty,3)$. This is easily done
numerically for e.g.~$C=0$. Notice next that by using a similar continuity
argument together with Lemma~\ref{l-3conj2} it is enough to verify condition
(ii) for $C=-1$ and for arbitrarily chosen values of $C$ in the intervals
$(-\infty,-1)$ and $(-1,3)$, respectively. This is again an easy numerical test
for e.g.~$C=-1$, $C=0$ and $C=2$. Finally, in order to check (iii) note that
at any point $z_0\in \Xi_{C}$ which is not a branching point of $\Gamma_C$
the roots of the asymptotic symbol equation~\eqref{aseq}
with maximal absolute value have necessarily distinct arguments since these
roots must be distinct. Fix any such $z_0$ and denote by $\tau_i(z_0)$,
$1\le i\le 3$, the (distinct) zeros of the corresponding asymptotic symbol
equation~\eqref{aseq}, so that $\tau_i(z_0)\neq 0$ for $1\le i\le 3$.
By~\eqref{eq:3conj} and~\eqref{eq:leadasymp} there exist complex constants
$K_i(z_0)$, $1\le i\le 3$, such that
$$p_n(z_0)=\sum_{i=1}^{3}K_i(z_0)\tau_i(z_0)^n\text{ for }n\in \{-2,-1,0\}
\cup \bN.$$
Using the initial values $p_{-2}(z)=p_{-1}(z)=0$ and $p_{0}(z)=1$ together with
the above formula for $p_n(z_0)$ one can easily check that $K_i(z_0)\neq 0$
for $1\le i\le 3$. Assume that $|\tau_1(z_0)|=|\tau_2(z_0)|\ge |\tau_3(z_0)|$.
Then
$$\frac{p_n(z_0)}{p_{n-1}(z_0)}=\tau_1(z_0)\frac{K_1(z_0)+K_2(z_0)
\left(\dfrac{\tau_2(z_0)}{\tau_1(z_0)}\right)^n+K_3(z_0)
\left(\dfrac{\tau_3(z_0)}{\tau_1(z_0)}\right)^n}{K_1(z_0)+K_2(z_0)
\left(\dfrac{\tau_2(z_0)}{\tau_1(z_0)}\right)^{n-1}+K_3(z_0)
\left(\dfrac{\tau_3(z_0)}{\tau_1(z_0)}\right)^{n-1}}$$
and since $\left|\dfrac{\tau_2(z_0)}{\tau_1(z_0)}\right|=1$ and
$\tau_2(z_0)\neq
\tau_1(z_0)$ it follows that the right-hand side of the above identity has no
limit as $n\rightarrow \infty$, which proves (iii).
\end{proof}

\section{Comparison of approximation schemes and open problems}

\subsection{Approximation via the algebraic scheme versus Pad\'e
approximation}\label{ss:adv}

As we already pointed out in the introduction, the usual Pad\'e approximation
and the approximation scheme introduced in this paper have an essentially 
different range of applications. Nevertheless, when applied to an arbitrarily
given algebraic function $f$ whose defining equation is known the scheme for 
approximating the branch with maximal modulus of $f$ proposed above has a 
number of advantages compared to the usual Pad\'e approximation scheme for $f$:

\begin{enumerate}
\item[(1)] While Pad\'e approximation requires the knowledge of the Taylor
expansion at $\infty$ of $f(z)$, our scheme uses
only the defining algebraic equation for $f(z)$.
\item[(2)] The regular poles of Pad\'e approximants concentrate on the union
of certain trajectories of the quadratic differential described in
\S 3. While it is in general rather difficult to understand the structure of
these
     trajectories, the induced equimodular discriminant
$\Xi_{f}$ has a more transparent
     definition and is also easier to study. Interesting examples of
     equimodular discriminants relevant to the behavior of Potts
     models on sequences of finite graphs tending to certain infinite
     lattices have been studied in a number of papers by A.~Sokal and his
     coauthors, see e.g.~\cite{So1} and~\cite{So2}. They also developed two
approaches for finding $\Xi_{f}$ explicitly. One is a more or less direct
     approach while the other is based on a resultant type
     equation which can also be found in~\cite{Bi}.
\item[(3)] In general, Pad\'e approximants have spurious poles with
uncontrolled behavior which makes uniform convergence impossible.
Our scheme has only a finite number of spurious poles tending to a
quite understandable finite set $\Si_{f,IN}$. In addition to that, the latter
scheme has also a well-controlled (exponential) rate of convergence.
\item[(4)] The denominators of the
rational approximants $\{r_{n}(z)\}_{n\in \bN}$
constructed in this paper satisfy a simple recursion
with  fixed rational coefficients.
By contrast, the denominators of Pad\'e approximants
satisfy a 3-term recurrence relation with varying coefficients which
are difficult to
calculate and whose behavior is apparently rather chaotic (see \S 5.2.3 below).
\end{enumerate}

\subsection{Related topics and open problems}

In this section we discuss
several questions pertaining to our approximation scheme that seem to have
important applications not only to the convergence theory of Pad\'e
approximants (\cite{St1}--\cite{St6}) but also to complex Sturm-Liouville
problems (\cite{BBS}), modern
statistical physics (\cite{So1,So2}) and the theory of (general)
orthogonal polynomials (\cite{Ch,Ne,ST,Sz}).

\subsubsection{}

Let $\{p_{n}(z)\}_{n\in \bN}$ be a family of
monic complex polynomials. In addition to the asymptotic ratio distribution
$\nu$ introduced in Definition~\ref{def3} one can associate to this family
an asymptotic probability measure $\mu$ defined as follows.

\begin{definition}\label{def11}
To a degree $d_{n}$ complex polynomial $p_{n}(z)$ we associate a finite
       probability measure $\mu_{n}$ called the {\em root counting measure}
       of $p_{n}(z)$ by placing the mass $\frac {1}{d_{n}}$ at any of its
simple roots.
If some root of $p_{n}(z)$ is multiple we place at this point
       a mass equal to the multiplicity of the root divided by $d_{n}$. The
limit
       $\lim_{n\to \infty}\mu_{n}$ (if is exists in the sense of the weak
convergence of measures) is called the {\em asymptotic root counting measure}
       of the family $\{p_{n}(z)\}_{n\in \bN}$ and denoted by $\mu$.
\end{definition}

It is clear from Definitions~\ref{def3} and~\ref{def11} that for any given
polynomial family $\{p_{n}(z)\}_{n\in \bN}$ the supports of $\nu$ and $\mu$
coincide.

\begin{problem}
Find the relation between the
asymptotic root counting measure $\mu$ and the asymptotic ratio distribution
$\nu$ for general polynomial families. Is it true that $\nu$ depends
only on $\mu$, i.e., can two polynomial families with the same
asymptotic root counting measure have different asymptotic ratio
distributions?
\end{problem}

A more concrete question about the
asymptotic root counting measure $\mu$ may be formulated as follows.
Take a sequence
$\big\{\big(Q_{1,n}(z),\ldots,Q_{k,n}(z)\big)\big\}_{n\in \bN}$
of $k$-tuples of
polynomials and consider the polynomial family $\{p_{n}(z)\}_{n\in \bN}$
satisfying the recurrence relation
$$p_{n+1}(z)=\sum_{i=1}^kQ_{i,n}(z)p_{n-i}(z)$$
with initial $k$-tuple
e.g.~$p_{-k+1}(z)=p_{-k+2}(z)=\ldots=p_{-1}(z)=0$, $p_{0}(z)=1$. We say that
the polynomial family $\{p_{n}(z)\}_{n\in \bN}$ belongs to the
{\em (generalized) Nevai class} if the sequence
$\big\{\big(Q_{1,n}(z),\ldots, Q_{k,n}(z)\big)\big\}_{n\in \bN}$
converges pointwise (and coefficientwise) to a fixed $k$-tuple of polynomials
$\big\{\big(\widetilde Q_{1}(z),\ldots,\widetilde Q_{k}(z)\big)\big\}$ when
$n\to\infty$.
Recall from~\cite{Ne} that the (ordinary) Nevai class of orthogonal
polynomials consists of all polynomial families $\{p_{n}(z)\}_{n\in \bN}$
satisfying a 3-term recursion of the form
$$zp_{n}(z)=a_{n+1}p_{n+1}(z)+b_{n}p_{n}(z)+a_{n}p_{n-1}(z),$$
where $a_{n}>0$, $b_{n}\in \bR$ for $n\in \bN$ and
$\lim_{n\to \infty}a_{n}=a\ge 0$, $\lim_{n\to\infty}b_{n}=b\in \bR$.
Theorem 5.2 of~\cite{Ne} claims that the asymptotic root counting measure
exists for any orthogonal polynomial family in Nevai's class and its density
is given by $\mu(x)=\frac{1}{\pi\sqrt{(\beta-x)(x-\alpha)}}$, where
$\alpha=b-2a,\;\beta=b+2a$. A classical result of Szeg\"o states
that in this case the density of $\nu$
     equals $\nu(x)=\frac {2\sqrt{(\beta-x)(x-\alpha)}}{\pi}$.

\begin{pb-prime}
Calculate the asymptotic root-counting measure $\mu$
(if it exists) for any polynomial family in the generalized Nevai class.
\end{pb-prime}

Problem $1'$ is important and nontrivial already in the context of
the present paper, namely for recurrences with fixed polynomial or
rational coefficients. Quite recently we were informed by A.~Sokal that
the asymptotic root-counting measure $\mu$ is conjecturally
given by the distributional derivative of $\log |y_{dom}(z)|$.

\subsubsection{}

Let us consider a finite recurrence relation with fixed polynomial
coefficients and arbitrarily prescribed initial values of the form
\begin{equation}\label{frec}
p_{n+1}(z)=\sum_{i=1}^kQ_{i}(z)p_{n-i}(z).
\end{equation}
Denote by $\Xi_Q$ the induced equimodular discriminant of the asymptotic
symbol equation of~\eqref{frec} viewed as a family of polynomials in the
variable $z$ (cf.~Definitions~\ref{def4} and~\ref{def7}). The following fact
was observed by the authors in extensive computer experiments and needs both a
mathematical ground and a conceptual explanation.

\begin{conjecture}
For any polynomial family $\{p_{n}(z)\}_{n\in \bN}$ satisfying
a finite recurrence relation of the form~\eqref{frec} such
that $\deg p_{n}(z)=n$ for $n\in \bN$
the zeros of any two consecutive polynomials
$p_{n+1}(z)$ and $p_{n}(z)$ interlace along the curve $\Xi_Q$ for all
sufficiently high degrees $n$.
\end{conjecture}

\begin{figure}[!htb]
\centerline{\hbox{\epsfysize=6.5cm\epsfbox{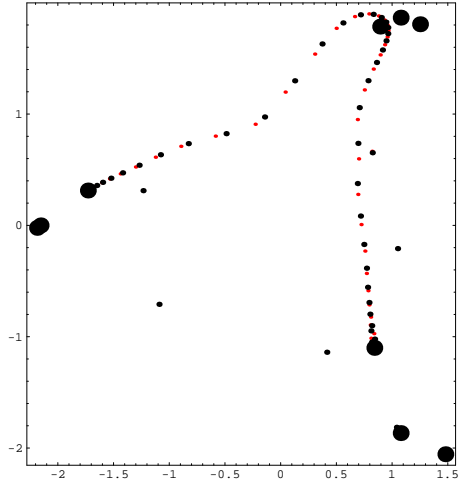}}}
\caption{Zeros of $p_{40}(z)$ and $p_{41}(z)$ satisfying the
4-term recurrence relation
$p_{n+1}(z)=(z+1-I)p_{n}(z)+(z+1)(z-I)p_{n-1}(z)+(z^3+10)p_{n-2}(z)$ with
$p_{0}(z)=z^6-z^4+I$, $p_{1}(z)=z-I+2$ and
$p_{2}(z)=(2+I)z^2-8$.}\label{fig3}
\end{figure}

An example illustrating the interlacing phenomenon conjectured above
is shown in Figure~\ref{fig3}. The large fat points in this picture are
as usual the branching points. Normal size fat points are the zeros of
$p_{41}(z)$ while thin points represent the zeros of $p_{40}(z)$. One
notices a typical interlacing pattern between these two sets of zeros.

A simple case of this phenomenon was studied in~\cite{ERR}
(see Conjecture~\ref{3Con}). Similar interlacing properties for complex zeros
along specified curves were noticed
in e.g.~\cite{BBS}. Some caution is actually required when defining the
interlacing property since the zeros of $p_{n}(z)$ do not lie
exactly on $\Xi_Q$. In the case when $\Xi_Q$ is a smooth curve
we may proceed as follows. Identify some sufficiently small neighborhood
$N(\Xi_Q)\subset \Omega$ of $\Xi_Q$ with the normal bundle
to $\Xi_Q$ by equipping $N(\Xi_Q)$ with a projection onto
$\Xi_Q$ along the fibers which are small curvilinear segments transversal
to $\Xi_Q$. Then one can say that two sets of points in
$N(\Xi_Q)$ {\em interlace} if their projections on $\Xi_Q$
interlace in the usual sense. If $\Xi_Q$ has singularities one
should first remove some sufficiently small neighborhoods of these
singularities and proceed in the above way on the remaining part of
$\Xi_Q$. The above conjecture then says that for any
sufficiently small neighborhood $N(\Xi_Q)$ of $\Xi_Q$ and any of
its identifications with the normal bundle to $\Xi_Q$ there exists
$n_{0}$ such that the interlacing property for
the roots of $p_{n}(z)$ and $p_{n+1}(z)$ holds for all $n\ge n_{0}$.
It is worth mentioning that the restriction $\deg p_{n}(z)=n$ is apparently
unnecessary.

\subsubsection{}

Consider the following analog of recursion~\eqref{eq:rec} with variable
coefficients
$$	-q_{n}(z)=\sum_{i=1}^{k}\frac{P_{k-i,n}(z)}{P_{k,n}(z)}q_{n-i}(z),$$
where $P_{j,n}(z)\in \bC[z]$, $0\le j\le k-1$, $n\in \bN$.
Assume that the sequence of $k$-tuples
$\{\big(P_{k,n}(z),P_{k-1,n}(z),\ldots,P_{0,n}(z)\big)\}_{n\in \bN}$ converges
(coefficientwise) uniformly on
compact subsets of $\bC$ to $\{\big(\widetilde
P_{k}(z),\widetilde P_{k-1}(z),\ldots, \widetilde P_{0}(z)\big)\}$. As in \S 1
we set $r_n(z)=\frac{q_n(z)}{q_{n-1}(z)}$ for  $n\in \bN$ and we denote by
$\widetilde y_{dom}(z)$ the dominant root of the limit equation
$$\sum_{i=0}^{k}\widetilde P_{k-i}(z)y^{k-i}=0.$$

\begin{problem}\label{pb2}
Is it true that the sequence of rational functions $\{r_{n}(z)\}_{n\in \bN}$
converges uniformly (in the complement of some set of zero Lebesgue
measure) to $\widetilde y_{dom}(z)$?
\end{problem}

Several results in this direction were obtained in \cite {BSh}. Note that an
affirmative answer to Problem~\ref{pb2} would imply in particular
that the coefficients of the
3-term recurrence relation~\eqref{eq:3term} can have no limits.

\subsubsection{}

Let $f(z)$ be an algebraic function and $\{r_{n}(z)\}_{n\in \bN}$ be a
sequence of rational approximants to $f(z)$ given by Theorem~\ref{th:Exist}.
It would be interesting to know what is the order of contact (in the sense
of~\eqref{eq:Pd}) of $r_{n}(z)$ and $f(z)$ at $\infty$. More generally,
in which cases does the sequence $\{r_{n}(z)\}_{n\in \bN}$ coincide (up to a
finite number of terms) with the sequence $\{[n/n](z)\}_{n\in \bN}$
of diagonal Pad\'e approximants to $f(z)$?

\subsubsection{}

A natural problem arising in the present context is the study of the
stratification of the space of algebraic functions according to the
topological structure of the induced equimodular discriminant $\Xi_{f}$. More
precisely, consider the space $\A_k$ of all algebraic
equations of degree $k$. An equation $P(y,z)\in \A_k$ is called {\em branching
nondegenerate} if all its branching points in the $z$-plane are simple.
As noted in \S 4 the endpoints of the curves in $\Xi_{f}$
can only occur at the branching points of $P(y,z)$ (cf.~\cite{Bi,BSh}). Thus
the topological structure of $\Xi_{f}$ can change if and only if either (i)
a pair of branching points collapses and $P(y,z)$ becomes branching
degenerate or (ii) a 4-tuple
of values of $f(z)$ have the same (maximal) absolute value. (Note
that for a given $z$ only the branching
points for which $|y(z)|$ is maximal are relevant to $\Xi_{f}$.)
What can be said about this stratification? In particular, how many
top-dimensional strata are there and what are their
(co)homology groups?
A similar
stratification related to certain quadratic differentials was studied in
\cite {Ba}.

\subsubsection{}

Well-known results on diagonal Pad\'e approximants at
$\infty$ going back to Markov \cite{Ma} claim that if the
analytic function under consideration is the Cauchy transform of a
positive measure supported on a real interval then the poles of all Pad\'e
approximants belong to this interval and thus both the Pad\'e Conjecture and
Nuttall's Conjecture are trivially valid for such functions (see also
\cite{Bu}). On the other hand, we saw in \S 3 that already algebraic
functions which are obtained as the Cauchy transform of a real but not
necessarily positive measure supported on an interval might violate
Nuttall's conjecture. It is therefore natural to investigate the Pad\'e
Conjecture and Nuttall's Conjecture for
algebraic functions representable as Cauchy transforms of positive
measures supported on compact subsets of $\bC$. An intriguing class of
such functions was recently discovered in \cite{BBSh}. Let $k\in \bN$ and
      consider the algebraic curve $\Ga$ given by the
          equation
          \begin{equation}
          \sum_{i=0}^k Q_{i}(z)w^i=0,
          \label{eq:basic}
          \end{equation}
where $Q_{i}(z)=\sum_{j=0}^{i}a_{i,j}z^j$ with $\deg Q_i\le i$ for
$0\le i\le k$. The curve $\Ga$ is called of
{\em general type} if the
          following two nondegeneracy requirements are satisfied:
\begin{enumerate}
\item[(i)] $\deg Q_{k}(z)=k;$
\item[(ii)] all roots of the (characteristic) equation
          \begin{equation}
          a_{k,k}+a_{k-1,k-1}t+\ldots+a_{0,0}t^k=0
          \label{eq:char}
          \end{equation}
have pairwise distinct arguments (in particular, $0$ is
          not a root of~\eqref{eq:char}).
\end{enumerate}

In \cite{BBSh} we showed that all the branches of an algebraic curve $\Ga$ of
general type vanish at $\infty$ and that each such branch is representable
(up to a constant factor) as the
	Cauchy transform of a certain compactly supported positive
	measure. We believe that the answer to the following question is
affirmative.

\begin{problem}
Is it true that all the branches of an arbitrary algebraic curve $\Ga$ of
general type satisfy both the Pad\'e
Conjecture and Nuttall's Conjecture near $\infty$?
\end{problem}

\subsection{The geometry of the algebraic scheme and
multidimensional continuous fractions}

Let us finally discuss certain similarities between the scheme for
approximating algebraic functions proposed in this paper and some recent
developments in multidimensional continuous
fractions, see \cite{Ar,Ko1,Ko2,La}. Recall first that the conjugates
of an algebraic number $\xi$ are by definition all algebraic
numbers other than $\xi$ satisfying the irreducible algebraic equation with
integer coefficients defining $\xi$. Notice
that Theorem~\ref{th:Exist} has an obvious analog for algebraic numbers.

\begin{definition}
An algebraic number $\xi$ is called {\em
dominant} if its absolute value is strictly larger than the absolute
values of all its conjugates.
\end{definition}

Clearly, any dominant algebraic number $\xi$ is necessarily real.
Assume additionally that  the
polynomial with integer coefficients
$P_{\xi}(t)=\al_{0}t^k+\al_{1}t^{k-1}+\ldots+\al_{k}$
defining $\xi$ is monic, that is, $\al_{0}=1$. (This assumption
is implicit in all the references quoted above where the defining polynomial is
always the characteristic polynomial of a matrix with integer
coefficients.) Consider the $k\times k$ matrix with integer
coefficients given by
$$M_{\xi}=\left(\begin{array}
{ccccc}-\al_{1}&-\al_{2}&-\al_{3}&\ldots &-\al_{k}\\
                      1     &0 &0      &\ldots & 0\\
                      0    &1 &0       & \ldots &0\\
		    \vdots& \vdots & \vdots &\vdots&\vdots\\
		    0    &  0     &\ldots &1&0
		    \end{array}\right). $$
   The characteristic polynomial of $M_{\xi}$ coincides with
   $P_{\xi}(t)$. Set $\bv_{0}=(1,0,\ldots,0)^t$ and let
   $\bv_{n}=\big(v_{n}^{(k)},v_{n}^{(k-1)},\ldots,v_{n}^{(1)}\big)^t
:=M_{\xi}^n\bv_{0}\in \bZ^k$ for
$n\in \bN$. Using the same methods as in \S 2 one can easily show that the
sequence of vectors $\{\bv_n\}_{n\in \bN}$ thus defined
enjoys the following property:

\begin{lemma}
For any $i\in \{1,\ldots,k-1\}$ one has
$\lim_{n\to\infty} \frac {v_{n}^{(i+1)}}{v_{n}^{(i)}}=\xi$.
\end{lemma}

Matrices with integer coefficients whose characteristic polynomials have only
distinct positive zeros are called {\em hyperbolic}. Such matrices are
the main ingredient in the vast program of studying multidimensional
continuous fractions suggested by V.~Arnold (following the ideas of
F.~Klein). More precisely, assuming that $M_{\xi}$ is a hyperbolic
matrix one gets a decomposition of $\bR^k$ into $2^k$ orthants spanned by all
  possible choices of directions given by the $k$
one-dimensional eigenspaces of $M_{\xi}$. The {\em sail} of such an orthant is
the convex hull of the set of all integer points in $\bR^k$ contained in
the orthant. If all the eigenspaces of $M_{\xi}$ are
irrational then Dirichlet's theorem on unities implies
that the  subgroup $G_{M_{\xi}}\subset GL_{k}(\bZ)$ acting on
$\bR^k$ and preserving all the eigenspaces of $M_{\xi}$ is isomorphic to
$\bZ_{2}\times \bZ^{k-1}$. The combinatorial properties of the
sails defined by $M_{\xi}$ and the action of the group $G_{M_{\xi}}$ on
these sails are the subject of
interesting and delicate generalizations of Lagrange's theorem on the
periodicity of continuous fractions for quadratic irrationalities. A natural
problem that arises in this context is as follows.

\begin{problem}
Let $\xi$ be a dominant algebraic number and assume that its associated
matrix $M_{\xi}$ is hyperbolic. Study the behavior of the sequence
$\{\bv_{n}\}_{n\in \bN}$ with respect to the union of the sails defined by
$M_{\xi}$.
\end{problem}


\begin{thebibliography}{99}

\bibitem{Ar}
V.~Arnold, Continuous fractions. (Russian), MCCME, 40pp, 2001.

\bibitem {BG-M}
G.~Baker, P.~Graves-Morris, {\em Pad\'e Approximants}, I and II,
Encyclopedia of Mathematics and Its Applications, Vols. 13 and 14,
Addison-Wesley, 1981.

\bibitem{Ba1} 
D.~Barrios Rolan\'{\i}a, G.~L\'opez Lagomasino,
{\em  Ratio asymptotics for polynomials orthogonal on arcs of the
unit circle}, Constr. Approx. {\bf 15} (1999), 1--31.

\bibitem{Ba2} 
D.~Barrios Rolan\'{\i}a, G.~L\'opez Lagomasino, E.~Torrano, 
{\em  Distribution of zeros and the asymptotics of polynomials
that satisfy three-term recurrent relations with complex coefficients}
(Russian), Mat. Sb. {\bf 184} (1993), 63--92;  translation in
Russian Acad. Sci. Sb. Math. {\bf 80} (1995), 309--333.

\bibitem{Ba}
Yu.~Baryshnikov, {\em On Stokes sets}, in
``New developments in singularity theory (Cambridge, 2000)'',
pp. 65--86, NATO Sci. Ser. II Math. Phys. Chem., Vol. 21, Kluwer Acad. Publ.,
Dordrecht, 2001.

\bibitem{BBS}
C.~M.~Bender, S.~Boettcher, V.~M.~Savage, {\em Conjecture on the interlacing
of zeros in complex Sturm-Liouville problems}, J. Math. Phys. {\bf 41} (2000),
6381--6387.

\bibitem{BKW1}
S.~Beraha, J.~Kahane, N.~J.~Weiss, {\em Limits of zeroes of recursively
defined polynomials}, Proc. Nat. Acad. USA {\bf 72} (1975), 4209.

\bibitem{BKW2}
S.~Beraha, J.~Kahane, N.~J.~Weiss, {\em Limits of zeros of recursively
defined families of polynomials}, in ``Studies in Foundations and
Combinatorics'', pp. 213--232, Advances in Mathematics Supplementary Studies
Vol. 1, ed. G.-C.~Rota, Academic Press, New York, 1978.

\bibitem{BeG}
C.~Berenstein, R.~Gay, Complex variables. An introduction. Grad. Texts in
Math., Vol. 125, Springer-Verlag, New York, 1991.

\bibitem{Bi}
N.~Biggs, {\em Equimodular curves}, Discrete Math. {\bf 259} (2002), 37--57.

\bibitem{BBSh}
J.~Borcea, R.~B\o gvad, B.~Shapiro, {\em Asymptotics of polynomial
eigenfunctions for linear ordinary differential operators, I: homogenized
spectral pencils}, manuscript in preparation.

\bibitem{BSh}
J.~Borcea, B.~Shapiro, {\em On the asymptotic ratio of recursively defined
families of functions}, manuscript in preparation.

\bibitem{Bu}
V.~I.~Buslaev, {\em On the Baker-Gammel-Wills conjecture in the theory of
Pad\'e approximants} (Russian), Mat. Sb. {\bf 193} (2002), 25--38;
translation in Sb. Math. {\bf 193} (2002), 811--823.

\bibitem{Bu2}
V.~I.~Buslaev, {\em On Poincar\'e's theorem and its applications to problems of
the convergence of continued fractions} (Russian),  Mat. Sb. {\bf 189} (1998),
13--28;  translation in  Sb. Math. {\bf 189} (1998),  1749--1764.

\bibitem{Bu3}
V.~I.~Buslaev, {\em Remarks on Poincar\'e and Van Vleck theorems},
Univ. of North Carolina Industr. Math. Inst. Res. Reports 99:14, (1999).

\bibitem{Ch}
T.~S.~Chihara, An introduction to orthogonal polynomials,
Gordon \& Breach, New York, 1978.

\bibitem{Du}
S.~Dumas, {\em Sur le d\'eveloppement des fonctions
elliptiques en fractions continues}, Th\`ese, Facult\'e de
Philosophie, Universit\'e de Zurich, 1908.

\bibitem{ERR}
\"O.~Egecioglu, T.~Redmond, C.~Ryavec, {\em From a polynomial Riemann
hypothesis to alternating sign matrices}, Electron. J. Combin.
{\bf 8} (2001), Research Paper 36, 51 pp (electronic).

\bibitem{Gel}
A.~O.~Gel'fond, Differenzenrechnung, Hochschulb\"ucher f\"ur Mathematik,
Vol. 41, VEB Deutscher Verlag der
Wissenschaften, Berlin, 1958.

\bibitem{Ka}
T.~Kato, Perturbation theory for linear operators, Die Grundlehren der
mathematischen Wissenschaften, Band 132, Springer-Verlag, New York, 1966.

\bibitem{Ko1}
E.~Korkina, {\em La p\'eriodicit\'e des fractions
continues multidimensionnelles}, C.R. Acad. Sci. Paris {\bf 319} (1994),
777-780.

\bibitem{Ko2}
E.~Korkina, {\em The simplest 2-dimensional continued
fraction}, J. Math. Sci. {\bf 82} (1996), 3680-3685.

\bibitem{La} G.~Lachaud, {\em Poly\`edre d'Arnol'd et voile d'un
c\^one simplicial: analogues du th\'eor\`eme de Lagrange}, C.R. Acad.
Sci. Paris {\bf 317} (1993), 711-716.

\bibitem{Lu}
D.~Lubinsky, {\em Rogers-Ramanujan and the Baker-Gammel-Wills (Pad\'e)
conjecture}, Ann. of Math. (2) {\bf 157} (2003), 847--889.

\bibitem{Ma}
A.~A.~Markoff [Markov], {\em Deux d\'emonstrations de la convergence de
certaines fractions continues}, Acta Math. {\bf 19} (1895), 93--104.

\bibitem{Ne}
P.~Nevai, Orthogonal polynomials, 4th edition, Mem.
Amer. Math. Soc., Vol. 213, 1979.

\bibitem{Ra}
E.~A.~Rakhmanov, {\em The asymptotic behavior of the ratio of orthogonal
polynomials} (Russian), I: Mat. Sb. {\bf 103(145)} (1977), 237-252; 
II: ibid. {\bf 118(160)} (1982), 104--117, 143. 

\bibitem{So1}
J.~Salas, A.~Sokal,
{\em Transfer matrices and partition-function zeros for
antiferromagnetic Potts models I}, J. Statist. Phys. {\bf 104} (2001),
609-699; http://arxiv.org/archive/cond-mat/0004330.

\bibitem{So2} A.~Sokal,
{\em Chromatic roots are dense in the whole complex plane}, Combin. Probab.
Comput.{\bf 13} (2004), 221-261;  http://arxiv.org/archive/cond-mat/0012369.

\bibitem{St1}
H.~Stahl, {\em The convergence of Pad\'e approximants to functions with
branch points}, J. Approx. Theory {\bf 91} (1997), 139-204.

\bibitem{St2}
H.~Stahl, {\em On the divergence of certain Pad\'e approximants and
the behavior of the associated orthogonal polynomials}, in:
``Polyn\^omes Orthogonaux et Applications'', pp. 321-330, Lect. Notes Math.,
Vol. 1171, Springer, Heidelberg, 1983.

\bibitem{St3}
H.~Stahl, {\em Divergence of diagonal Pad\'e approximants and the
asymptotic behavior of orthogonal polynomials associated with
nonpositive measures}, Constr. Approx.{\bf 1} (1985), 249-270.

\bibitem{St4}
H.~Stahl, {\em Diagonal Pad\'e approximants to hyperelliptic functions},
Ann. Fac. Sci. Toulouse, n$^o$ sp\'ecial Stieltjes, pp. 121-193, 1996.

\bibitem{St5}
H.~Stahl, {\em Spurious poles in Pad\'e approximation}, J. Comput. Appl.
Math.{\bf 99} (1998), 511-527.

\bibitem {St6}
H.~Stahl, {\em Extremal domains associated with an analytic function,
I-III}, Complex Var. Theory Appl. {\bf 4} (1985), 311-354.

\bibitem {ST}
H.~Stahl, W.~Totik, General orthogonal polynomials,
Encyclopedia of Mathematics, Vol. 43, Cambridge Univ. Press,
Cambridge, 1992.

\bibitem{St}
R.~Stanley, Enumerative Combinatorics, Vol. 1, Wadsworth
Brooks/Cole
Math. Ser., The Wadsworth \& Brooks/Cole Advanced Books \& Software,
Monterey, CA, 1986.

\bibitem{Sz}
G.~Szeg\"o, Orthogonal Polynomials, 4th edition, Coll. Publ., Vol. XXIII,
Amer. Math. Soc., Providence, RI, 1975.

\end{thebibliography}
\end{document}